\documentstyle[amssymb,amsfonts]{amsart}

\newcommand{\Ecal}{{\cal E}}

\def\be#1{\begin{equation} \label{#1}}
\def\bs{\begin{split}}
\def\bi{\begin{itemize}}
\def\es{\end{split}}
\def\ba{\begin{align}}
\def\bas{\begin{align*}}
\def\ea{\end{align}}
\def\eas{\end{align*}}

\def\R{{\mbox{\bf R}}}
\def\T{{\mbox{\bf T}}}
\def\Tcal{{\mbox{\bf T}}}
\def\m{m}
\def\dist{{\mbox{\rm dist}}}
\def\diam{{\mbox{\rm diam}}}
\def\supp{{\mbox{\rm supp}}}
\def\margin{{\mbox{\rm margin}}}

\def\Qcal{{\cal Q}}
\def\K{{2^{k(\frac{1}{p}-\frac{1}{2})}}}
\def\D{{\mbox{\bf D}}}
\def\Z{{\mbox{\bf Z}}}
\def\eps{\varepsilon}

\def\emph#1{{\it #1}}
\def\textbf#1{{\bf #1}}
\newenvironment{proof}{\noindent {\bf Proof} }{\endprf\par}
\def \endprf{\hfill  {\vrule height6pt width6pt depth0pt}\medskip}

\theoremstyle{plain}
  \newtheorem{theorem}[subsection]{Theorem}
  \newtheorem{proposition}[subsection]{Proposition}
  \newtheorem{lemma}[subsection]{Lemma}
  \newtheorem{corollary}[subsection]{Corollary}
  \newtheorem{conjecture}[subsection]{Conjecture}
  \newtheorem{problem}[subsection]{Problem}

\theoremstyle{remark}

\theoremstyle{definition}
  \newtheorem{definition}[subsection]{Definition}

\include{psfig}

\begin{document}

\title[Endpoint restriction for the cone]{Endpoint
bilinear restriction theorems for the cone, and some sharp null form estimates}

\author{Terence Tao}
\address{Department of Mathematics, UCLA, Los Angeles, CA 90024}
\email{tao@@math.ucla.edu}

\subjclass{42B15, 35L05}
\begin{abstract}
Recently Wolff \cite{wolff:cone} obtained a nearly sharp $L^2$ bilinear
restriction theorem for bounded subsets of the cone in general dimension.  We obtain the endpoint of Wolff's estimate and generalize to the case when one of the subsets is large.  As a consequence, we are able to deduce some nearly-sharp
$L^p$ null form estimates.
\end{abstract}

\maketitle

\section{Introduction}

Let $n \geq 2$ be a fixed integer.  We say that a function $\phi: \R^{n+1} \to H$ is a \emph{red wave} if it takes values in a finite dimensional complex Hilbert space $H$, and its space-time Fourier transform
$\hat \phi$ is an $L^2$ measure on the set
$$ 2^k \Sigma^{red} := \{ (\xi,|\xi|): \angle(\xi,e_1) \leq \frac{\pi}{8}, 2^k \leq |\xi| \leq 2^{k+1} \}$$
for some integer $k$, where $e_1$ is a fixed basis vector. 
Similarly, we say that $\psi: \R^{n+1} \to H'$ is a \emph{blue wave} if
it takes values in a finite dimensional complex Hilbert space $H'$ and
$\hat \psi$ is an $L^2$ measure on
$$ 2^k \Sigma^{blue} := \{ (\xi,-|\xi|): \angle(\xi,e_1) \leq \frac{\pi}{8},
2^k \leq |\xi| \leq 2^{k+1} \}$$
for some integer $k$.
In both cases we call $2^k$ the \emph{frequency} of the waves $\phi$, $\psi$.

Red and blue waves both solve the free wave equation, but propagate along different sets of characteristics.  Note that blue waves are the time reversal of red waves.  Also, these waves are automatically smooth and bounded thanks to the frequency localization.  The vector valued formulation will be convenient for inductive reasons, but our implicit constants shall always be independent of $H$ and $H'$.

We define the \emph{energy} of $\phi$, $\psi$ by
\be{energy-def} E(\phi) := \| \phi(t) \|_2^2, \quad E(\psi) := \| \psi(t)\|_2^2
\end{equation}
where $t \in \R$ is arbitrary.  This definition is independent of the choice of $t$, and is related to the standard notion of energy by the formula
$$ 2^{2k} E(\phi) \sim \frac{1}{2} \int_{\R^n} |\phi_t(0)|^2 + |\nabla \phi(0)|^2\ dx.$$

Throughout the paper, $p_0 = p_0(n)$ will denote the exponent
\be{p0-def}
p_0 := \frac{n+3}{n+1}.
\end{equation}
The main result of this paper is the following bilinear estimate.  

\begin{theorem}\label{main}  Let $\phi$ be a red wave of frequency 1, and
$\psi$ be a blue wave of frequency $2^k$ for some $k \geq 0$. Then we have
\be{est} \| \phi \psi \|_p \lesssim 2^{k(\frac{1}{p}-\frac{1}{2}+\eps)} E(\phi)^{1/2} E(\psi)^{1/2}
\end{equation}
for all $2 \geq p \geq p_0$, $\eps > 0$.  In particular, $\phi$, $\psi$ have frequency 1 then
\be{est-0}
\| \phi \psi \|_p \lesssim E(\phi)^{1/2} E(\psi)^{1/2}.
\end{equation}
\end{theorem}

In the above theorem and in the sequel, the implicit constants may depend on $\eps$ but are independent of $H$, $H'$, and $\phi \psi: \R^{n+1} \to H \otimes H'$ denotes the tensor product $\phi \otimes \psi$ of $\phi$ and $\psi$.

The estimate \eqref{est-0} solves a conjecture of Machedon and Klainerman.   The restriction $p \geq p_0$ is sharp; see e.g.
\cite{damiano:null}, \cite{tv:cone1}, \cite{tvv:bilinear}.  For $2 \leq p \leq \infty$ the theory is much simpler, and the best possible estimate is
$$ \| \phi \psi \|_p \lesssim 2^{kn(\frac{1}{2}-\frac{1}{p})} E(\phi)^{1/2} E(\psi)^{1/2}.$$
This is easily proved from the $p=2$ case and Sobolev embedding.

The estimate \eqref{est-0} is a genuinely bilinear estimate and cannot be proven directly from linear estimates.  Indeed,
the Strichartz estimate \cite{strichartz:restrictionquadratic}
combined with the H\"older inequality only yields the range 
$p \geq (n+1)/(n-1)$, while Plancherel's theorem and Cauchy-Schwarz
only gives the range $p \geq 2$ (see e.g. \cite{kl-mac:null}).  In the $n=2$ case, the fact that
one could go below $p=2$ was first shown by Bourgain \cite{borg:cone},
and in \cite{tv:cone1} a concrete range of $p$ was given, namely
$p > 2 - \frac{8}{121}$.  More recently, Wolff \cite{wolff:cone} obtained the range $p > p_0$ for all dimensions $n$.  Thus \eqref{est-0} is the endpoint of that in \cite{wolff:cone}.

The generalization \eqref{est} of \eqref{est-0} is necessary in order to develop sharp null form estimates, as we shall see in Section \ref{final}.  The estimate \eqref{est} is sharp except for the $\eps$.  To see this, let
$\psi$ be a blue wave whose Fourier transform is supported in a unit ``square''
in $2^k \Sigma_{blue}$, and which is comparable to 1 on a $2^k \times 1$ tube oriented in a blue null direction, and let $\phi = (\phi_i)_{i=1}^{2^k}$ be a vector-valued
red wave of frequency $1$ such that $\phi_i$ is comparable to 1 on $B_i$, where the $B_i$ are a family of $2^k$ unit balls that cover the above tube.   The $\eps$ term can probably be removed when $p > p_0$.  When $p=p_0$ the author conjectures that $\eps$ can still be removed, but this seems to require an extremely delicate analysis and some new Kakeya estimates for null rays.  However, it should be fairly straightforward to replace the $2^{\eps k}$ factor by a polynomial $k^C$ in this case. 

Broadly, our strategy to prove Theorem \ref{main} is as follows.  We shall localize the estimate \eqref{est} to a cube $Q$ of side-length $R \gg 2^k$, and obtain a bound independent of $R$.  This will be obtained by induction on $R$, as follows.  

If $\phi$ and $\psi$ are dispersed fairly evenly throughout $Q$, we shall decompose $Q$ into sub-cubes of side-length $2^{-C_0} R$ and decompose $\phi$ and $\psi$ into smaller waves, each of which is concentrated on one of these sub-cubes.  By an interpolation between bilinear $L^1$ and $L^2$ estimates as in Wolff \cite{wolff:cone}, the cross-terms are well controlled, and one can replace $\phi$ and $\psi$ by a ``quilt'' of waves on the $2^{-C_0} R$-cubes.  One then applies the induction hypothesis to each sub-cube and sums up. 

This tactic works well when $\phi$ and $\psi$ are dispersed, but there is a problem when almost all the energy of $\phi$ and $\psi$ simultaneously concentrate in a disk $D$ of radius $r \ll R$.  By Huygens' principle the wave $\phi\psi$ is then concentrated in the double light cone generated by $D$.  Restricting to this smaller set, we can exploit a more favourable bilinear $L^1$ estimate (Corollary \ref{doublecone}) than the more trivial bilinear $L^1$ estimate \eqref{l1-noloc} used in the non-concentrated case.  One can then repeat the non-concentrated argument, and localize $\phi$ and $\psi$ all the way down to the scale of $r$, at which point the waves become non-concentrated and one obtains enough of a gain to close the induction.

The proof is unfortunately rather complex.  In an attempt to give the reader a sense of the full argument without drowning in the details, the author has abstracted the argument into several sections.  We first give the top-level argument, in which Theorem \ref{main} is deduced from several major propositions.  Then, we give the medium-level argument, in which these major propositions are deduced from some elementary estimates and a key estimate, Proposition \ref{lung}.  Finally, we devote several sections to the proof of Proposition \ref{lung}.  We have tried to make each section as self-contained as possible, so that the arguments in each section rely only on the Propositions and Lemmata of the previous sections, and not on the method of arguments or on notation specific to a single section.

The author wishes to thank Damiano Foschi, Sergiu Klainerman, Christoph Thiele, and Tom Wolff for many helpful discussions, and for making their preprints readily available.  This work was conducted at UCLA and at UNSW, with additional support from NSF grant DMS-9706764 and grants from the Sloan and Packard foundations.  The author is also indebted to the referee for finding some errors in the first version of this paper (which claimed \eqref{est} without the $2^{\eps k}$ factor).

\section{Some notation}\label{notation}

In this section we list some notation which will be used throughout the argument.

We fix $N \gg 1$ to be a large integer depending only on $n$
($N = 2^{n^{10}}$ will suffice); the disclaimer ``assuming $N$ is sufficiently large'' will be implicit throughout our arguments.  We also let $0 < \eps \ll 1$ be an arbitrary small number. We also let $C_0$ denote an integer much larger than $N$, $1/\eps$; for instance,
\be{c0-def}
C_0 := 2^{\lfloor N/\eps \rfloor^{10}}
\end{equation}
will suffice.  Generally speaking, we use $N$ as a large exponent, and $C_0$ as a very large constant (large enough to dominate any reasonable quantity arising from $N$ or $\eps$).  As a rule of thumb, any term containing an $r^{-N}$ or $R^{-N}$ factor may be ignored for all practical purposes; these error terms only arise because one cannot quite simultaneously localize in both space and frequency.

We use the notation $A := B$ to indicate that $A$ is being defined to equal $B$.  We shall extend this notation in several ways, for instance $\hat f := g$ indicates that $f$ is being defined via the Fourier transform.

We let $C$ denote various large numbers that vary from line to line (possibly depending on $N$, $\eps$, but will not depend explicitly on $C_0$), and let $A \lesssim B$ or $A=O(B)$ denote the estimate $A \leq CB$ where $C$ depends only on $n$ and $\eps$.  Similarly we use $A \ll B$ to denote $A \leq C^{-1} B$.  Generally speaking, we shall use the $\lesssim$ notation to control error terms, but will need the more precise $\leq$ notation for the main terms due to the inductive nature of the argument.  In particular, we will be dealing with many estimates of the form $A \leq B + CE$, $A \leq (1+Cc)B + c^{-C}E$, or $A \leq (1+CNc)B + c^{-C} E$, where $B$ is the main term, $E$ is the error term, and $0 < c \ll 1$ is some small parameter which we may optimize in later.  For such estimates it will be very important that the factor in front of $B$ is either 1 or very close to 1, as we will be unable to close the induction otherwise.

Note that our constants $C$ are independent of the dimension of the spaces $H$, $H'$.  This will be important for the induction argument.

If $\phi(x,t)$ is a function of both space and time, we use $\phi(t)$ to denote the spatial function $\phi(t)(x) := \phi(x,t)$.

We use the hat notation for the spatial Fourier transform
$$ \hat f(\xi) := \int_{\R^n} e^{-2\pi i x \cdot \xi} f(x)\ dx$$
as well as the spacetime Fourier transform
$$ \hat \phi(\xi,\tau) := \int_{\R^{n+1}} e^{-2\pi i (x \cdot \xi + t \tau)} \phi(x,t)\ dx dt,$$
with the meaning being clear from context.  We define the \emph{frequency support} $\supp(\hat f)$ of a function $f$ to be the support of the Fourier transform $\hat f$.

We shall always treat $\R^{n+1}$ as endowed with the Euclidean metric and
never with the Minkowski metric, so that terms such as length $|x|$,
angle $\angle(x,y)$, etc. retain their usual meaning in $\R^{n+1}$.  On the other hand, we will employ the Lorentz transforms on occasion, especially when we derive null form estimates in Section \ref{final}.

We now define some familiar geometric objects, namely disks, cubes, cones, and conic neighbourhoods.

A \emph{disk} will be any subset $D$ of $\R^{n+1}$ of the form
$$ D = D(x_D,t_D;r_D) := \{ (x,t_D): |x-x_D| \leq r_D \}$$
for some $(x_D,t_D) \in \R^{n+1}$ and $r_D > 0$.  The reader should note that disks are $n$-dimensional objects even though they reside in $\R^{n+1}$.  We call $t_D$ the \emph{time co-ordinate} of the disk $D$.  If $D$ is a disk, we define the cutoff function $\tilde \chi_D$ on $\R^{n+1}$ by
\be{tchib-def} \tilde \chi_D(x,t) := (1 + \frac{|x-x_D|}{r_D})^{-N^{10}}.
\end{equation}
If $D = D(x_D,t_D;r_D)$ is a disk and $c > 0$, we define $cD$ to be the disk
$cD = D(x_D,t_D;cr_D)$, and the \emph{disk exterior} $D^{ext} = D^{ext}(x_D,t_D;r_D)$ to be the region
$$ D^{ext}(x_D,t_D;r_D) := \{ (x,t_D): |x-x_D| > r_D \}.$$
We endow disks and disk exteriors with spatial Lebesgue measure $dx$.

We define $Q(x_Q,t_Q;r_Q)$ to be the $n+1$-dimensional cube in $\R^{n+1}$ centered at $(x_Q,t_Q)$ with side-length $r_Q$ and with sides parallel to the axes.  We call the interval $[t_Q - r_Q/2, t_Q + r_Q/2]$ the \emph{lifespan} of $Q$.
If $Q = Q(x_Q,t_Q;r_Q)$ is a cube and $c>0$, we use $cQ$ to denote the cube
$cQ := Q(x_Q,t_Q;cr_Q)$.  Finally, we define the cubical annuli $Q^{ann}(x_Q,t_Q;r_1,r_2)$ by
$$ Q^{ann}(x_Q,t_Q;r_1,r_2) := Q(x_Q,t_Q;r_2) \backslash Q(x_Q,t_Q;r_1).$$

Let $\underline \Sigma = \underline \Sigma_n$ denote the spatial region
$$ \underline \Sigma := \{ \xi: \frac{1}{2} \leq |\xi| \leq 4, \angle(\xi, e_1) \leq \frac{\pi}{4} \}.$$
If $(x_0,t_0) \in \R^{n+1}$, define the \emph{red cone} with vertex $(x_0,t_0)$
to be the set
$$ C^{red}(x_0,t_0) := \{ (x_0 + r\omega, t_0 - r): r \in \R, \omega \in S^{n-1} \cap \underline\Sigma \},$$
and the \emph{blue cone} at this vertex to be
$$ C^{blue}(x_0,t_0) := \{ (x_0 + r\omega, t_0 + r): r \in \R, \omega \in S^{n-1} \cap \underline\Sigma \}.$$
For any $r \gg 1$, we define $C^{red}(x_0,t_0;r)$, $C^{blue}(x_0,t_0;r)$ to be the $r$-neighbourhoods of $C^{red}(x_0,t_0)$, $C^{blue}(x_0,t_0)$ respectively.  
Finally, we define the combined conic neighbourhood $C^{purple}(x_0,t_0;r)$ to be 
$$ C^{purple}(x_0,t_0;r) := C^{red}(x_0,t_0;r) \cup C^{blue}(x_0,t_0;r).$$

For any integer $j$, let $\D_j$ denote the dilation operator
$$ \D_j \phi(x,t) := \phi(2^j x, 2^j t),$$
and $\T$ to be the time reversal operator
$$ \T \phi(x,t) := \phi(x,-t).$$
The operator $\T$ maps red waves onto blue waves and vice versa, while the operator $\D_j$ maps waves of frequency $1$ onto waves of frequency $2^j$.

In our induction argument we shall frequently be decomposing red and blue waves into smaller waves.  Unfortunately, these decompositions often enlarge the Fourier support of the waves slightly, which is a potential obstruction to closing the induction.  To get around this we introduce the notion of the \emph{margin} $\margin(\phi)$ of a wave $\phi$.  More precisely, if $\phi$ is a red wave of frequency $1$, we define $\margin(\phi)$ to be the non-negative real number
$$ \margin(\phi) := \dist(\supp(\hat \phi), \Sigma^+ \backslash \Sigma^{red})$$
where $\Sigma^\pm$ denotes the light cones
$$ \Sigma^\pm := \{ (\xi, \pm |\xi|): \xi \in \R^n \}.$$
We extend this definition to general red or blue waves by the formulae
$$ \margin(\T \D_k \phi) := \margin(\D_k \phi) := \margin(\phi).$$
Thus, for instance, if $\psi$ is a blue wave of frequency $2^k$, then 
$$ \margin(\psi) := 2^{-k} \dist(\supp(\hat \psi), 2^k (\Sigma^- \backslash \Sigma^{blue})).$$
The concept of margin is only needed to overcome the technical obstruction mentioned earlier, and otherwise plays no role of importance.  

\section{The top-level proof of Theorem \ref{main}}\label{top-sec}

In this section we give the top-level proof of Theorem \ref{main}, in which we state the key propositions and give the inductive argument which yields the Theorem from these propositions.  We shall prove these propositions in later sections.

Throughout the proof $n \geq 2$, $k \geq 0$, $\eps > 0$ will be fixed.
It suffices to verify the case $k \geq 2^{C_0^{10}}$, since the $k < 2^{C_0^{10}}$ case follows by applying a Lorentz transform.  We shall thus assume $k \geq 2^{C_0^{10}}$ throughout the argument.  In particular, any factor depending on $k$ will dwarf any factor depending on $C_0$, which in turn dwarfs any factor depending on $N$, which in turn dwarfs any absolute constants.

We can also assume that $k$ is a multiple of $C_0$.  This is a purely notational convenience that we shall only use in Section \ref{fundamental}.

It is known (see e.g.
\cite{borg:cone}, \cite{kl-mac:null}, \cite{tv:cone1}, \cite{mock:cone}, \cite{damiano:null}, or Lemma \ref{mock} below) that 
\be{l2-noloc}
\| \phi \psi \|_2 \lesssim E(\phi)^{1/2} E(\psi)^{1/2}.
\end{equation}
Thus Theorem \ref{main} holds for $p=2$.  By interpolation it thus suffices to prove the theorem when 
\be{p-def}
p = p_0 = \frac{n+3}{n+1},
\end{equation}
and we shall implicitly assume this throughout the proof. 

\begin{definition}\label{ar0-Def}  For any $R \geq C_0 2^k$, we define $A(R)$ to be the best constant for which the inequality
\be{ar-def} \| \phi \psi \|_{L^{p}(Q_R)} \leq A(R) E(\phi)^{1/2} E(\psi)^{1/2}\end{equation}
holds for all spacetime cubes $Q_R$ of side-length $R$, red waves $\phi$ of frequency 1, and blue waves $\psi$ of frequency $2^k$ such that one has the margin requirement
\be{margin-bound}
\margin(\phi), \margin(\psi) \geq 1/100 - (2^k/R)^{1/N}.
\end{equation}
The waves $\phi$, $\psi$ may take values in arbitrary finite-dimensional Hilbert spaces.
\end{definition}

Note that it suffices to verify \eqref{ar-def} for those $\phi$, $\psi$ which obey the normalization
\be{energy-bound}
E(\phi)=E(\psi)=1.
\end{equation}

From \eqref{l2-noloc} it is clear that $A(R)$ is finite for each $R$.
The margin requirements on $A(R)$ are technicalities which are needed for the induction on $R$ to work properly, as many of our decompositions will decrease the margin of $\phi$ and $\psi$ slightly.  However, we may remove the margin requirements by a finite partition of space and frequency, and some mild Lorentz transforms to obtain the crude estimate
\be{crude-eq}
\| \phi \psi \|_{Q_R} \lesssim A(R') E(\phi)^{1/2} E(\psi)^{1/2}
\end{equation}
for any cube $Q_R$ of side-length $R$, any $R' \sim R$, and any red and blue waves $\phi$, $\psi$ of frequency $1$, $2^k$ respectively.  In particular to prove \eqref{est} it suffices to show that
\be{ark}
A(R) \lesssim 2^{CC_0} 2^{\eps k} \K
\end{equation}
uniformly for all $R \geq C_0 2^k$.

Because of the increasingly strict margin requirements as $R \to \infty$
we see that $A(R)$ is not necessarily increasing in $R$.  We therefore define the auxiliary quantity $\overline A(R)$ for all $R \geq C_0 2^k$ by
$$ \overline A(R) := \sup_{C_0 2^k \leq r \leq R} A(r).$$

We now begin the proof of \eqref{ark}.  In Section \ref{initial-proof} we shall prove \eqref{ark} when $R$ is close to $2^k$:

\begin{proposition}\label{initial}
For any $R \geq C_0 2^k$, we have the bound
$$A(R) \lesssim 2^{CC_0} (R/2^k)^C \K.$$
\end{proposition}

This proposition is needed to begin the inductive argument.  The exact power of $(R/2^k)$ is unimportant as we shall soon improve this bound substantially.

In Section \ref{upshift-proof}, we adapt the localization ideas from \cite{wolff:cone} to prove the following recursive inequality on $A(R)$.

\begin{proposition}\label{upshift}
Let $R \geq 2C_0 2^k$ and $0 < c \leq 2^{-C_0}$, and let $\phi$, $\psi$ be re and blue waves of frequency 1 and $2^k$ respectively which obey the relaxed margin requirement
\be{margin-relaxed}
\margin(\phi), \margin(\psi) \geq 1/100 - 2(2^k/R)^{1/N}.
\end{equation} 
Then for any cube $Q_R$ of side-length $R$ we have
\be{upshift-bound}
\| \phi \psi \|_{L^p(Q_R)} \leq [(1+Cc) \overline A(R/2) + c^{-C} \K] E(\phi)^{1/2} E(\psi)^{1/2}.
\end{equation}
In particular, we have
$$ A(R) \leq (1+CNc) \overline A(R/2) + c^{-C} \K.$$
\end{proposition}

The slight relaxation of margin requirements in Proposition \ref{upshift} (as compared to \eqref{margin-bound}) is important for technical reasons, but should be ignored for a first reading.

The loss of $CNc$ in the main term comes about because we need to divide an $R$-cube into sub-cubes at various scales, and then shrink those cubes further by about $c$ in order to create some separation between those sub-cubes (this is necessary, otherwise the error terms can blow up).  In practice we can always optimize $c$ so that its impact on the argument is negligible.

By setting $c = (2^k/R)^{\frac{\eps}{CN}}$ in Proposition \ref{upshift}, we obtain
$$ A(R) \leq (1+CN(2^k/R)^{\frac{\eps}{CN}}) \overline A(R/2) + C (2^k/R)^{\eps/N} \K$$
for all $R \geq 2C_0 2^k$.  Iterating this and using Proposition \ref{initial} when $R \sim 2^k$, we obtain
\be{a-grow}
A(R) \lesssim 2^{CC_0} (R/2^k)^{\eps/N} \K
\end{equation} 
for any $R \geq C_0 2^k$.  This can then be used to obtain Theorem \ref{main} for $p > p_0$ (cf. \cite{wolff:cone}, \cite{tv:cone2}, \cite{borg:cone}), but we shall not do so here.  

From the preceding discussion we observe that \eqref{ark} is already proven for $R \leq 2^{C_0 k}$.  Thus we may assume that $R > 2^{C_0 k}$.  $R$ is now very large, dominating most quantities which depend only on $k$, $C_0$, $N$, and $\eps$.

In the large $R$ case we need to introduce the notion of energy concentration.  

\begin{definition}\label{conc-def}
For any $r > 0$, spacetime cube $Q$, red wave $\phi$, and blue wave $\psi$, we define the \emph{energy concentration} $E_{r,Q}(\phi,\psi)$ to be the quantity
$$ E_{r,Q}(\phi,\psi) := \max(\frac{1}{2} E(\phi)^{1/2} E(\psi)^{1/2}, \sup_D \| \phi \|_{L^2(D)}\| \psi \|_{L^2(D)})$$
where $D$ ranges over all disks of radius $r$ whose time co-ordinate $t_D$ is contained in the lifespan of $Q$.
\end{definition}

Clearly we have $E_{r,Q}(\phi,\psi) \sim E(\phi)^{1/2} E(\psi)^{1/2}$ for any $r$,
with equality or near-equality only occurring when $\phi$, $\psi$ simultaneously concentrate in a disk of radius $r$.  However, the choice of whether to use 
$E_{r,Q}(\phi,\psi)$ instead of $E(\phi)^{1/2} E(\psi)^{1/2}$ will be crucial to make a certain induction work.

Roughly speaking, the strategy for the $R \geq 2^{C_0 k}$ case is to show that Proposition \ref{upshift} can be improved slightly unless there is substantial energy concentration.  This slight improvement will be enough to close the induction, however, one must still deal with the concentrated case.  In this case we use Huygens principle to restrict ourselves to a small region of a double cone $C^{purple}$, in which case one can obtain an improvement of Proposition \ref{upshift} by other means.

To make this strategy more precise we shall need the following technical variant of $A(R)$ which is weighted slightly to exploit the gain in the non-concentrated case.

\begin{definition}\label{ar-Def}  For any $R \geq 2^{C_0 k}$ and any $r', r > 0$, we define $A(R,r,r')$ to be the best constant for which the inequality
\be{arr-def} \| \phi \psi \|_{L^{p}(Q_R \cap C^{purple}(x_0,t_0; r'))} \leq A(R,r,r') (E(\phi)^{1/2} E(\psi)^{1/2})^{1/p} E_{r,C_0 Q_R}(\phi,\psi)^{1/p'}
\end{equation}
holds for all spacetime cubes $Q_R$ of side-length $R$, all $(x_0,t_0) \in \R^{n+1}$, red waves $\phi$ of frequency 1, and blue waves $\psi$ of frequency $2^k$ such that one has the margin requirement \eqref{margin-bound}
and the energy bound \eqref{energy-bound}.
\end{definition}

As before, it suffices to verify \eqref{arr-def} assuming the energy normalization \eqref{energy-bound}.  It is important that the right-hand side if \eqref{arr-def} is exactly as stated, and not (for instance) the comparable quantity $A(R,r,r') E(\phi)^{1/2} E(\psi)^{1/2}$. In practice $r$ and $r'$ shall usually be comparable in size.

The preceding heuristics regarding concentration can be formalized in the following Proposition, which we prove in Section \ref{back-sec}.  This Proposition is basically an application of Huygens' principle, combined with some more sophisticated arguments to deal with the highly concentrated case $r \lesssim R^{1/2 + 4/N}$.

\begin{proposition}\label{back-prop}
For any $R \geq 2^{C_0 k}$, we have
$$ A(R) \leq (1 - C_0^{-C}) \sup_{2^{C_0 k} \leq \tilde R \leq R; \tilde R^{1/2+4/N} \leq r} A(\tilde R,r,C_0(1+r)) + 2^{CC_0} 2^{\eps k} \K.$$
\end{proposition}

The requirement $\tilde R^{1/2 + 4/N} \leq r$ is somewhat difficult to obtain, but it is necessary to do so because the tools we shall develop to control $A(R,r,r')$ have a spatial uncertainty of about $\sqrt{R}$ and so one cannot effectively exploit any concentration effects near or below this scale. This uncertainty of $\sqrt{R}$ is responsible for all the powers of $r^{1/N}$ and $R^{1/N}$ which appear in the arguments; these powers should be ignored for a first reading.  The key point to observe in Proposition \ref{back-prop} is that we have somehow wrested a small gain $(1-C_0^{-C})$ from the main term on the right-hand side, thanks to the beneficial effects of non-concentration in \eqref{arr-def}.  This gain will allow us to absorb all error terms and $(1+Cc)$ factors in the other Propositions in this section, thus closing the induction.

To use this inequality inductively we need to bound $A(R,r,r')$ in terms of $A(R)$.  
From \eqref{crude-eq} it is easy to show that $A(R,r,r') \leq C A(R)$, but this is too crude to close the induction, and one must take some more care with the constant in the leading term.

The bounding of $A(R,r,r')$ by $A(R)$ can be split into two stages.  First, we shall use finite speed of propagation in Section \ref{flip-sec} to observe that one can obtain the desired bound in the non-concentrated case $r \geq C_0^C R$:

\begin{proposition}\label{r-flip}
For any $R \geq 2^{C_0 k}$, $r \geq C_0^C R$, $r' > 0$, and $0 < c \leq 2^{-C_0}$, we have
$$ A(R,r,r') \leq (1 + Cc) \overline{A(R)} + c^{-C} \K.$$
\end{proposition}

For the concentrated case we iterate the following Proposition, which is proven in \ref{medium-sec}.

\begin{proposition}\label{r-medium}
For any $R \geq 2^{C_0 k}$ and $C_0^C R \geq r > R^{1/2 + 3/N}$, we have
$$ A(R,r,r') \leq (1 + CNc) A(R/C_0,r(1 - C r^{-1/3N}),r') + c^{-C} (1 + \frac{R}{2^k r})^{-1/N} \K$$
for any $0 < c \leq 2^{-C_0}$.
\end{proposition}

The decay of $(1 + \frac{R}{2^k r})^{-1/N}$ in the error term is crucial for this endpoint result as it allows us to avoid losing a logarithmic factor $\log(R/r)$ in the iteration, which would otherwise be fatal to the proof of the endpoint.  This decay ultimately arises from the improved energy estimates on cones as encoded in Lemma \ref{bluecone}.  However, the presence of the $2^k$ means that we still lose a factor of $k$ when summing over scales $R$, and this is the main source of the $2^{\eps k}$ loss in \eqref{est}.

These two Propositions combine to give

\begin{corollary}\label{r-total}
For any $R \geq 2^{C_0 k}$ and $r \geq R^{1/2+4/N}$, we have
$$ A(R,r,C_0(1+r)) \leq (1 + CNc) \overline{A(R)} + c^{-C} 2^{\eps k} \K$$
for any $0 < c \leq 2^{-C_0}$.
\end{corollary}

\begin{proof}
We may assume that $r < C_0^C R$ since the claim follows from Proposition \ref{r-flip} otherwise.  Let $J$ be the first integer such that $r \geq 2^{-J} C_0^C R$; from the hypotheses we have $J \lesssim \log(r)$.  Define $r =: r_0 >  r_1 > \ldots > r_J$ inductively by $r_{j+1} := r_j(1 - Cr_j^{1-1/3N})$.  One can verify inductively that $r_j = r + O(j r^{-1/4N})$ for all $j$, and in particular that $r_j \sim r$.

From Proposition \ref{r-medium} we have
\begin{align*}
 A(R/2^j,r_j,C_0(1+r)) &\leq (1 + CNc_j) A(R/2^{j+1}, r_{j+1}, C_0(1+r))\\
&+ c_j^{-C} (1 + \frac{R}{2^j 2^k r})^{-1/N} \K
\end{align*}
for any $0 < c_j \leq 2^{-C_0}$.  

Observe that
$$ 1 + \frac{R}{2^j 2^k r} \gtrsim C_0^{-C} 2^{(J-j-k)_+}.$$
If we thus set $c_j := k^{-1} c 2^{-(J-j-k)_+/(2CN)}$ for a suitable constant $C$, we thus obtain
\begin{align*} A(R/2^j,r_j,C_0(1+r)) \leq 
&(1 + CN k^{-1} c 2^{-(J-j-k)_+/(2CN)}) A(R/2^{j+1}, r_{j+1}, C_0(1+r))\\
& + k^C c^{-C} 2^{-(J-j-k)_+/2N} \K.
\end{align*}
Iterating this, one obtains
$$ A(R,r,C_0(1+r)) \leq (1 + CNc) A(R/2^J, r_J, C_0(1+r)) + k^C c^{-C} \K.$$
The claim then follows from Proposition \ref{r-flip}.
\end{proof}

Combining this with \eqref{back-prop} and setting $c := 2^{-C_0}$, we thus obtain
$$ A(R) \leq (1 - C_0^{-C}) \overline{A(R)} + 2^{CC_0} 2^{\eps k} \K$$
for all $R \geq 2^{C_0 k}$.  Combining this with \eqref{a-grow} we thus see that this inequality thus holds for all $R \geq C_0 2^k$.  Taking suprema and using the monotonicity of $\overline{A(R)}$ we thus obtain
$$ \overline{A(R)} \leq (1 - C_0^{-C}) \overline{A(R)} + 2^{CC_0} 2^{\eps k} \K$$
for all $R \geq C_0 2^k$, which implies
$$ \overline{A(R)} \leq 2^{CC_0} 2^{\eps k} \K,$$
and \eqref{ark} follows for all $R \geq C_0 2^k$ as desired.

It remains to prove Propositions \ref{initial}-\ref{r-medium}.  It turns out that most of these propositions follow as consequences of a localization property for red and blue waves, Proposition \ref{lung}, which we shall state in the next section, after some notation.  In Sections \ref{initial-proof}-\ref{flip-sec} we show how this Proposition implies Propositions \ref{initial}-\ref{r-medium}.  Finally, in Sections \ref{opposite-sec}-\ref{fundamental} we give a proof of Proposition \ref{lung}.

\section{The main proposition}

In this section we state the main proposition of the argument, Proposition \ref{lung}.  As we shall see in the next few sections, this Proposition will be the main tool used to prove Propositions \ref{initial}, \ref{upshift}, and \ref{r-medium}, and also plays a minor role in the proof of Proposition \ref{r-flip}.

Proposition \ref{lung} involves the localization of waves $\phi$, $\psi$ on a large cube $Q$ into smaller waves localized on sub-cubes of $Q$.  To make this precise we must introduce some notation.

If $Q$ is a cube of side-length $R$, and $j \geq 0$ is an integer, we may partition $Q$ into $2^{(n+1)j}$ cubes of side-length $2^{-j} R$; we use $\Qcal_j(Q)$ to denote the collection of these cubes.

If $Q$ is a cube and $j \geq 0$ is an integer, we define a \emph{red wave table
$\phi$ on $Q$ with depth $j$} to be any red wave with the vector form
$$ \phi =: (\phi^{(q)})_{q \in \Qcal_j(Q)},$$
where the components $\phi^{(q)}$ may themselves be vector-valued.  
If $0 \leq j' < j$, and $q' \in \Qcal_{j'}(Q)$, we define $\phi^{(q')}$
to be the red wave table on $q'$ with depth $j-j'$ given by
$$ (\phi^{(q')})^{(q)} := \phi^{(q)} \hbox{ for all } q \in \Qcal_{j-j'}(q').$$
Note that
$$ E(\phi) = \sum_{q' \in \Qcal_{j'}(Q)} E(\phi^{(q')})$$
for all $0 \leq j' \leq j$.

If $\phi$ is a red wave table on $Q$ of depth $j$ and $0 \leq j' \leq j$, we define the \emph{$j'$-quilt} $[\phi]_{j'}$ of $\phi$ to be the non-negative function
$$ [\phi]_{j'} := \sum_{q \in \Qcal_{j'}(Q)} |\phi^{(q)}|\ \chi_q.$$
Note we have the pointwise estimates
\be{quilt-est}
|\phi^{(q)}| \chi_q \leq [\phi]_j \leq [\phi]_{j-1} 
\ldots \leq [\phi]_0 = |\phi| \chi_Q
\end{equation}
for all $q \in \Qcal_j(Q)$.
We define blue wave tables and their quilts analogously.  

The estimates \eqref{quilt-est} are of course very crude, and we shall frequently be exploiting various improvements to this estimate in the sequel.

If $Q$ is a cube, $k \geq 0$ is an integer, and $0 \leq c \ll 1$, we define 
the \emph{$(c,k)$-interior} $I^{c,k}(Q)$ of $Q$ by
\be{ick-def}
I^{c,k}(Q) := \bigcup_{q \in \Qcal_k(Q)} (1-c)q.
\end{equation}
The advantage of working with $I^{c,k}(Q)$ instead of $Q$ is that the sub-cubes of $I^{c,k}(Q)$ have some significant separation properties.

With these notational preliminaries, we are now ready to state Proposition \ref{lung}.

\begin{proposition}\label{lung}  
Let $R \geq C_0 2^k$, $0 < c \leq 2^{-C_0}$, and let $\phi$, $\psi$ be red and blue waves with frequency 1 and $2^k$ respectively, which obey the energy normalization \eqref{energy-bound}
and the relaxed margin requirement \eqref{margin-relaxed}.  
For any cube $Q$, define the set $X(Q) \subset Q$ by
\be{x-def}
X(Q) := \bigcap_{j=C_0}^k I^{c 2^{-(k-j)/N}, j}(Q).
\end{equation}
Then for any cube $Q$ of side-length $C R$, we can find a red wave table $\Phi$ on $Q$ with depth $k$ and frequency 1, and a blue wave table $\Psi$ on $Q$ with depth $C_0$ and frequency $2^k$, such that we have the margin estimates
\be{marge-est}
\margin(\Phi), \margin(\Psi) \geq 1/100 - (2^{k+C_0}/R)^{1/N},
\end{equation}
 we have the energy estimates
\be{energy-est}
E(\Phi), E(\Psi) \leq 1+Cc,
\end{equation}
and we have the inequality
\be{triad}
 \| \phi \psi \|_{L^p(X(Q))} \leq \| [\Phi]_k [\Psi]_{C_0} \|_{L^p(X(Q))}
+ c^{-C} \K.
\end{equation}
For any cone $C^{purple}(x_0,t_0;r)$ with $r > 1$, we may improve \eqref{triad} to
\be{triad-cone}
\begin{split}
\| \phi \psi \|_{L^p(X(Q) \cap C^{purple}(x_0,t_0;r))} \leq &\| [\Phi]_k [\Psi]_{C_0} \|_{L^p(X(Q) \cap C^{purple}(x_0,t_0;r))}\\
&+ c^{-C} \K (1 + \frac{R}{2^k r})^{-1/N}.
\end{split}
\end{equation}
Furthermore, we have the persistence of non-concentration
\be{energy-persist}
E_{r(1 - C_0 r^{-1/4N}),2Q}(\Phi, \Psi) \leq E_{r,2Q}(\phi,\psi) + Cc + c^{-C} R^{C-N/2}
\end{equation}
for all $r \gtrsim R^{1/2 + 3/N}$. 
\end{proposition}

The reason why we use $X(Q)$ instead of $Q$ is that the sub-cubes of $X(Q)$ have some non-zero separation between them, which will let us avoid some unpleasantness in Section \ref{packet}.  To return from $X(Q)$ to $Q$ we shall use an averaging lemma, Lemma \ref{averaging-lemma}; this causes the loss of $Cc$ in many of the Propositions stated previously.  The estimates \eqref{energy-est}, \eqref{energy-persist} are stating that $\Phi$, $\Psi$ are ``smaller than or equal to'' $\phi$, $\psi$ in energy norm, while \eqref{triad}, \eqref{triad-cone} state that the quilts $[\Phi]_k$, $[\Psi]_{C_0}$ are good approximations for $\phi$, $\psi$ respectively.  This replacement of waves by quilts of essentially equal or lesser energy allows one to induct efficiently (providing that one has constants close to 1 in the main terms).

Note that the margin requirements on $\phi$, $\psi$ are slightly weaker than those in Definition \ref{ar-Def}.  The gain of $(1 + \frac{R}{2^k r})^{-1/N}$ in \eqref{triad-cone} over \eqref{triad} is responsible for the corresponding gain in Proposition \ref{r-medium} as compared against similar estimates such as Proposition \ref{upshift}.  This gain is essential as it needs to compensate both for the $Cc$ loss in \eqref{energy-est} and for the $\log(R)$ loss that would otherwise arise from an induction on scale.

Interestingly, the low-frequency wave $\phi$ can be localized much further than the high-frequency wave $\psi$; the former can be localized to cubes of $2^{-k}$ the side-length, whereas the high-frequency waves can only be localized by about $2^{-C_0}$ before the error estimates begin to deteriorate logarithmically.  This behaviour is also seen in the example given in the introduction demonstrating that the power of $k$ in \eqref{est} is essentially sharp.  In applications we shall only need to localize to $2^{-C_0}$, however in the proof of Proposition \ref{lung} we shall need to localize $\phi$ all the way down to $2^{-k}$ before one can begin to localize $\psi$.

The proof of Proposition \ref{lung} is the longest part of the argument, and is basically a ``pigeonhole-free'' version of the arguments in Wolff \cite{wolff:cone}.  Since it is the statement of this Proposition, rather than the method of proof, which are required to prove the Propositions of the previous section, we shall defer the proof of Proposition \ref{lung} to sections \ref{bilinear-sec}-\ref{fundamental}.  For now, we devote ourselves to the question of how this Proposition, combined with some other more elementary tools, can be used to prove Propositions \ref{initial}-\ref{r-medium}.

\section{Energy estimates}

In this section we record some fairly easy energy estimates which will be used throughout the paper.

Let $\phi$, $\psi$ be red and blue waves.
By integrating \eqref{energy-def} along the life-span of a cube we see that
\be{L2-triv} 
\| \phi \|_{L^2(Q_R)} \lesssim R^{1/2} E(\phi)^{1/2}, \quad
\| \psi \|_{L^2(Q_R)} \lesssim R^{1/2} E(\psi)^{1/2}
\end{equation}
for any cube $Q_R$ of side-length $R$.  By H\"older's inequality we thus have
\be{l1-noloc} 
\| \phi \psi \|_{L^1(Q_R)} \lesssim R E(\phi)^{1/2} E(\psi)^{1/2}.
\end{equation}
To compensate for the loss of $R$ in \eqref{l1-noloc} we shall often seek variants of \eqref{l2-noloc} in which one has a gain of $R^{-(n-1)/4}$.  This gain exactly balances the loss of $R$ when interpolated at the exponent $p=p_0$.  One such example of this gain is Lemma \ref{mock}; see also \eqref{non-conc2}, \eqref{l2-cloc}.

To tackle the case $k \gg 1$ of widely differing frequencies we shall need the following improvement of \eqref{l1-noloc} when $\phi$ is replaced by a quilt. 

\begin{lemma}\label{phipsi}  Let $Q_R$ be a cube of side-length $R > 0$, $j \geq 0$ is an integer, 
$\phi$ be a red wave table on $Q_R$ with depth at least $j$, and $\psi$ be a blue wave. Then
$$ \|[\phi]_j \|_2 \lesssim 2^{-j/2} R^{1/2} E(\phi)^{1/2},$$
and
$$
\| [\phi]_j \psi \|_1 \lesssim 2^{-j/2} R E(\phi)^{1/2} E(\psi)^{1/2}.
$$
\end{lemma}

\begin{proof}  Let $q \in \Qcal_j(Q_R)$.  From \eqref{L2-triv} for $\phi^{(q)}$ we have
$$ \|\phi^{(q)}\|_{L^2(q)} \lesssim 2^{-j/2} R^{1/2} E(\phi^{(q)})^{1/2}.$$
Square summing this in $q$ we obtain the first estimate,
and the second estimate follows from \eqref{L2-triv} for $\psi$ and H\"older's inequality.
\end{proof}

This improvement over \eqref{l1-noloc} is one demonstration of the gain involved when passing from a wave to a quilt.  The other major advantage of quilts is that they allow one to localize estimates on large cubes to estimates on small cubes so efficiently that one can make induction on scale arguments work.

\section{An averaging lemma}

The following averaging lemma will be needed in the proofs of Proposition \ref{initial}, \ref{upshift}, \ref{r-medium}.  It allows one to average out a bounds on $X(Q)$ to obtain a bound on $Q$.

\begin{lemma}\label{averaging-lemma}
Let $R > 0$, $0 < c \leq 2^{-C_0}$, and $F$ be an arbitrary smooth function.  Let $Q_R$ be a cube of side-length $R$.
Then there exists a cube $Q$ of side-length $CR$ contained in $C^2 Q_R$ such that
$$ \| F \|_{L^p(Q_R)} \leq (1 + CNc) \| F \|_{L^p(X(Q)))},$$
where $X(Q)$ is the set in \eqref{x-def}.
\end{lemma}

\begin{proof}
By the pigeonhole principle it suffices to show that
$$ \| F \|_{L^p(Q_R)}^p \leq \frac{1}{|Q_R|}
\int_{Q_R} (1 + Cc)^p \| F \|_{L^p(Q_R \cap (X(Q(x_0,t_0;CR))))}^p\ dx_0 dt_0.$$
From the symmetry and translation covariance of $X(Q)$ we have the identity
$$
\int_{Q_R} \| F \|_{L^p(Q_R \cap (X(Q(x_0,t_0;CR))))}^p\ dx_0 dt_0
= \int_{Q_R} |F(x,t)|^p |X(Q(x,t;CR)) \cap Q_R|\ dx dt.$$
On the other hand, from \eqref{ick-def} we have
$$
|Q(x_0,t_0;CR) \backslash I^{c,k}(Q(x_0,t_0;CR))| \lesssim c |Q(x_0,t_0;CR)|
$$
so from \eqref{x-def} we have
$$
|Q(x,t;CR) \backslash X(Q(x,t;CR))| \lesssim Nc |Q(x,t;CR)|$$
and thus that
$$ |Q_R| \leq (1 + CNc) |X(Q(x,t;CR)) \cap Q_R|.$$
The claim then follows.
\end{proof}

\section{Proof of Proposition \ref{initial}}\label{initial-proof}

We now have enough machinery to prove Propositions \ref{initial}, \ref{upshift}, and \ref{r-medium}.  We begin with Proposition \ref{initial}.

Fix $R \geq C_0 2^k$, and let $\phi$, $\psi$ be red and blue waves of frequency 1 and $2^k$ respectively, obeying the margin requirement \eqref{margin-bound} and the energy normalization \eqref{energy-bound}.  Let $Q_R$ be a cube of side-length $R$.  To prove Proposition \ref{initial} it suffices to show that
$$ \| \phi \psi \|_{L^p(Q_R)} \lesssim 2^{CC_0} (R/2^k)^C \K.$$

Set $c := 2^{-C_0}$, so that $1 + CNc \sim 1$.
From Lemma \ref{averaging-lemma} with $F := \phi \psi$, we may find a cube $Q$ of side-length $CR$ such that
$$
\| \phi \psi \|_{L^p(Q_R)} \lesssim \| \phi \psi \|_{L^p(X(Q))}.$$
Let $\Phi$, $\Psi$ be as in Proposition \ref{lung}.
By \eqref{triad} and \eqref{quilt-est}, we thus have
$$
\| \phi \psi \|_{L^p(Q_R)} \lesssim \| [\Phi]_k \Psi \|_{p}
+ 2^{C C_0} \K.$$
By \eqref{quilt-est}, \eqref{l2-noloc} and \eqref{energy-est} we have
$$ \| [\Phi]_k \Psi \|_{2} \leq \| \Phi \Psi \|_2 
\lesssim E(\Phi)^{1/2} E(\Psi)^{1/2} \lesssim 1,$$
while from Lemma \ref{phipsi} and \eqref{energy-est} we have
$$ \| [\Phi]_k \Psi \|_{1} \leq 2^{-k/2} R E(\Phi)^{1/2} E(\Psi)^{1/2}
\lesssim 2^{-k/2} R = (R/2^k) 2^{k/2}.$$
Interpolating between the two, one obtains
$$\| [\Phi]_k \Psi \|_p \lesssim (R/2^k)^C \K,$$
and the claim follows.
\endprf

\section{Proof of Proposition \ref{upshift}}\label{upshift-proof}

Let $\phi$, $\psi$, $Q_R$, $c$ be as in Proposition \ref{upshift}.  We may assume the normalization \eqref{margin-bound}.
We may assume that $R \geq 2^{CC_0} 2^k$ for some large $C$, since the claim \eqref{upshift-bound} follows from Proposition \ref{initial} otherwise.

From Lemma \ref{averaging-lemma} with $F := \phi \psi$, we may find a cube $Q$ of side-length $CR$ such that
$$
\| \phi \psi \|_{L^p(Q_R)} \leq (1+ CNc) \| \phi \psi \|_{L^p(X(Q))}.$$
Let $\Phi$, $\Psi$ be as in Proposition \ref{lung}.  By \eqref{triad}, we thus have
$$
\| \phi \psi \|_{L^p(Q_R)} \leq (1+ CNc) \| [\Phi]_k [\Psi]_{C_0} \|_{L^p(X(Q))}
+ c^{-C} \K.$$
By \eqref{quilt-est} and the hypothesis $k \gg C_0$, we thus have
\be{skiff}
\| \phi \psi \|_{L^p(Q_R)} \leq (1+ CNc) \| [\Phi]_{C_0} [\Psi]_{C_0} \|_{L^p(Q)}
+ c^{-C} \K.
\end{equation}
We expand
$$ \| [\Phi]_{C_0} [\Psi]_{C_0} \|_{L^p(Q)} = (\sum_{q \in \Qcal_{C_0}(Q)}
\| \Phi^{(q)} \Psi^{(q)} \|_{L^p(q)}^p)^{1/p}.$$
By \eqref{ar-def} and \eqref{marge-est}, then the inclusion $l^p \subset l^1$, followed by Cauchy-Schwarz and then \eqref{energy-est} we thus have
\begin{align*}
\| [\Phi]_{C_0} [\Psi]_{C_0} \|_{L^p(Q)} &\leq A(2^{-C_0}R)
(\sum_{q \in \Qcal_{C_0}(Q)}
E(\Phi^{(q)})^{p/2} E(\Psi^{(q)})^{p/2})^{1/p} \\
&\leq \overline A(R/2)
\sum_{q \in \Qcal_{C_0}(Q)}
E(\Phi^{(q)})^{1/2} E(\Psi^{(q)})^{1/2} \\
&\leq \overline A(R/2) 
(\sum_{q \in \Qcal_{C_0}(Q)} E(\Phi^{(q)}))^{1/2}
(\sum_{q \in \Qcal_{C_0}(Q)} E(\Psi^{(q)}))^{1/2}\\
&= \overline A(R/2) E(\Phi)^{1/2} E(\Psi)^{1/2}\\
&\leq (1+Cc) \overline A(R/2).
\end{align*}
Inserting this back into \eqref{skiff} we obtain the desired estimate \eqref{upshift-bound}.
\endprf

\section{Proof of Proposition \ref{r-medium}}\label{medium-sec}

We now prove Proposition \ref{r-medium}.  This will be a reprise of the argument used to prove Proposition \ref{upshift}, but will use non-concentration to improve upon the inclusion $l^p \subset l^1$ used previously.

Fix $R \geq 2^{C_0 k}$, $0 < c \leq 2^{-C_0}$, $r' > 0$, and $C_0^C R \geq r > R^{1/2 + 3/N}$.  Let $\phi$, $\psi$ be red and blue waves of frequency 1, $2^k$ respectively obeying \eqref{margin-relaxed}.  We may assume the normalization \eqref{energy-bound}.  Let $Q_R$ be a spacetime cube of side-length $R$, and let $(x_0,t_0)$ be a point in spacetime.  To prove Proposition \ref{r-medium} it suffices to show that
\begin{align*}
\| \phi \psi \|_{L^{p}(Q_R \cap C^{purple}(x_0,t_0; r'))} 
\leq& (1 + CNc) A(R/C_0,\tilde r,r')
E_{r,C_0 Q_R}(\phi,\psi)^{1/p}\\
& + c^{-C} (1 + \frac{R}{2^k r'})^{-1/N} \K
\end{align*}
where $\tilde r := r(1 - C_0 r^{-1/3N})$.
Here we have used the fact that $E_{r,C_0 Q_R}(\phi,\psi)^{1/p}
 \sim 1$.

Applying Lemma \ref{averaging-lemma} with 
$$ F := \phi \psi \chi_{ C^{purple}(x_0,t_0; r') }$$
we see that there exists a cube $Q$ of side-length $CR$ contained in $C^2 Q_R$ such that
$$
\| \phi \psi \|_{L^{p}(Q_R \cap C^{purple}(x_0,t_0; r'))} 
\leq (1+CNc)
\| \phi \psi \|_{L^{p}(X(Q) \cap C^{purple}(x_0,t_0; r'))}.$$
Let $\Phi$, $\Psi$ be as in Proposition \ref{lung}.  By \eqref{triad-cone} and the inclusion $2Q \subset C_0 Q_R$, it thus suffices to show that
\begin{align*}
\| [\Phi]_k [\Psi]_{C_0} \|_{L^p(X(Q) \cap C^{purple}(x_0,t_0;r'))}
\leq &(1 + Cc) A(R/C_0,\tilde r,r')
E_{r,2Q}(\phi,\psi)^{1/p'} \\
&+ c^{-C} R^{C-N/2},
\end{align*}
since the latter factor can easily be absorbed into $c^{-C} (1 + \frac{R}{2^k r'})^{-1/N} \K$ since $R > 2^{C_0 k}$.

By \eqref{quilt-est} we can estimate $[\Phi]_k$ by $[\Phi]_{C_0}$.  By
\eqref{energy-persist} it thus suffices to show that
$$
\| [\Phi]_{C_0} [\Psi]_{C_0} \|_{L^p(Q \cap C^{purple}(x_0,t_0;r'))}
\leq (1+Cc) A(R/C_0,\tilde r,r') E_{\tilde r,2Q}(\Phi,\Psi)^{1/p'}.$$
Here we have used some trivial bound in $A(R/C_0,\tilde r, r')$ which is polynomial in $R$, e.g. by \eqref{l2-noloc} and H\"older.

Raising both sides to the $p^{th}$ power and expanding, we reduce to showing that
\be{medium-targ}
\sum_{q \in \Qcal_{C_0}(Q)}
\| \Phi^{(q)} \Psi^{(q)} \|_{L^p(q)}^p
\leq (1+Cc) A(R/C_0, \tilde r,r')^p E_{\tilde r,2Q}(\Phi,\Psi)^{p/p'}.
\end{equation}
The cube $q$ has side-length $2^{-C_0} C R \ll R/C_0$, so in particular $C_0 q \subset 2Q$.  By \eqref{arr-def} and \eqref{marge-est} we thus have 
$$ \| \Phi^{(q)} \Psi^{(q)} \|_{L^p(q)}
A(R/C_0, \tilde r, r') (E(\Phi^{(q)})^{1/2} E(\Psi^{(q)})^{1/2})^{1/p}
E_{\tilde r,2Q}(\Phi^{(q)},\Psi^{(q)})^{1/p'}.$$
Also, by \eqref{energy-est} we have
$$
E_{\tilde r,2Q}(\Phi^{(q)},\Psi^{(q)})
\leq (1 + Cc)
E_{\tilde r,2Q}(\Phi,\Psi).$$
Comparing these estimates to \eqref{medium-targ}, it suffices to show that
$$
\sum_{q \in \Qcal_{C_0}(Q)} E(\Phi^{(q)})^{1/2} E(\Psi^{(q)})^{1/2}
\leq 1 + Cc.$$
But this follows from Cauchy-Schwarz and \eqref{energy-est} (cf. the proof of Proposition \ref{upshift}).
\endprf

\section{Spatial localization, and Huygens' principle}\label{localize}

We shall shortly begin the proof of Propositions \ref{back-prop} and \ref{r-flip}, but we must first develop some machinery concerning the localization of a wave to a disk, and the resulting estimates arising from Huygens' principle.  This machinery will also be useful in proving the persistence of non-concentration estimate \eqref{conc-persist} in Section \ref{packet}.

To localize a wave to a disk we shall need a bump function $\eta_0$ and an evolution operator $U(t)$, which we now construct.

Let $\eta_0$ to denote a fixed non-negative Schwarz function on $\R^n$ with total mass 1 and whose Fourier transform is supported on the unit disk; such a function can be constructed for instance by $\eta_0 = \hat \varphi * \hat \varphi$, where $\varphi$ is a real even bump function supported near the origin such that $\varphi(0)=1$.  For any $r > 0$ we define $\eta_r$ by
\be{etar-def}
\eta_r(x) := r^{-n} \eta_0(\frac{x}{r}).
\end{equation}
The function $\eta_r$ is thus an $L^1$-normalized Schwarz function concentrating on an $r$-disk with good frequency localization properties.  We shall use this cutoff function both in this section and in the construction of the wave packet decomposition in Appendix I.

Fix $a(\xi)$ to be a bump function supported on $\underline \Sigma$ which equals 1 on the spatial projection of $\Sigma^{red}$ and $\Sigma^{blue}$.
We define the evolution operator $U(t) = U^{red}(t)$ by
\be{u-def}
\widehat{U(t)f}(\xi) := a(\xi) e^{2\pi i t|\xi|} \hat f(\xi)
\end{equation}
and the kernel $K_t$ by
$$ K_t(x) := \int a(\xi) e^{2\pi i (x \cdot \xi + t|\xi|)}\ d\xi.$$
We have the propagation law
\be{prop}
 \phi(t) = U(t) \phi(0) = \phi(0) * K_t
\end{equation}
for all red waves $\phi$ with frequency $1$, and all times $t_1, t_2$.  
A stationary phase computation gives the pointwise estimate
\be{k-bound} |K_t(x)| \lesssim (1 + \dist((x,t),C^{red}(0,0))^{-N^{10}}.
\end{equation}
We shall also use the operator $U(t)$ in the proof of an energy estimate on cones (Lemma \ref{bluecone}), and in the proof of wave packet decomposition (Lemma \ref{tube-decomp}) in Appendix I.

\begin{definition}  If $D = D(x_D,t_D;r)$ is a disk, and $\phi$ is a red wave of frequency 1, we define $P_D \phi$ at time $t_D$ by
$$ P_D \phi(t_D) := (\chi_D * \eta_{r^{1-1/N}}) \phi(t_D)$$
and at other times $t$ by
$$ P_D \phi(t) := U(t-t_D) P_D \phi(t_D).$$
Here $\eta_r$.
For red waves of frequency $2^j$ we define $P_D$ by the formula
$$ P_D \D_j \phi := \D_j P_D \phi,$$
and for blue waves we define $P_D$ by the formula
$$ P_D \T \phi := \T P_D \phi.$$
\end{definition}

The operator $P_D$ localizes a wave to the disk $D$ at time $t_D$, while $1-P_D$ similarly localizes to $D^{ext}$.  More precisely:

\begin{lemma}\label{cutoff}  Let $j$ be an integer, $r \geq C_0 2^{-j}$.  Let $D = D(x_D,t_D;r)$ be a disk with radius $r$, and let $\phi$ be a red wave with frequency $2^j$ and margin 
$\margin(\phi) \geq C_0 (2^j r)^{-1+1/N}$.  Then $P_D \phi$ is a red wave of frequency $2^j$ which satisfies margin estimate
$$ \margin(P_D \phi) \geq \margin(\phi) - C_0(2^j r)^{-1+1/N}$$
and the energy estimates 
\ba
\| \tilde \chi_D^{-N} P_D \phi \|_{L^2(D_+^{ext})} &\lesssim (2^j r)^{-N} E(\phi)^{1/2}
\label{essential-concentration}\\
\| (1-P_D) \phi \|_{L^2(D_-)} &\lesssim (2^j r)^{-N} E(\phi)^{1/2}
\label{essential-vanish}\\
E(P_D \phi) &\leq \| \phi \|_{L^2(D_+)}^2 + C (2^j r)^{-N} E(\phi)\label{local-energy}\\
E((1-P_D) \phi) &\leq \| \phi \|_{L^2(D_-^{ext})}^2 + C (2^j r)^{-N} E(\phi) \label{non-local-energy}\\
E(P_D \phi), E((1-P_D) \phi) &\leq E(\phi) \label{energy-minor}
\end{align}
where $D_-$, $D_+$ are the disks
$$ D_\pm := D(x_D,t_D; r(1 \pm (2^j r)^{-1/2N})).$$
\end{lemma}

Note in particular that
$$ \frac{1}{2}D \subset D_- \subset D \subset D_+ \subset 2D.$$

\begin{proof}  
By scaling it suffices to verify this when $j=0$.  The claims follow directly from the easily verified estimates
$$ 0 \leq \chi_D * \eta_{r^{1-1/N}}(x) \leq 1,$$
$$ \chi_D * \eta_{r^{1-1/N}}(x) \lesssim r^{-N} \hbox{ for } x \in D_-,$$
$$ 1 - C r^{-N} \leq \chi_D * \eta_{r^{1-1/N}}(x) \hbox{ for } x \in D_+^{ext}.$$
\end{proof}

Applying $\T$ we see that similar statements hold for blue waves.

We now investigate the localization properties of $P_D$ for times other than $t_D$.  Heuristically speaking, $P_D \phi$ is supported on the cone $C^{red}(x_D, t_D; Cr)$, while $(1-P_D) \phi$ vanishes on the cube $Q(x_D, t_D; C^{-1} r)$.  More precisely, we have

\begin{lemma}\label{huygens}  Let $D$ be a disk of radius $r \geq C_0$, and
let $\phi$ be a red wave of frequency $1$ and margin $\margin(\phi) \geq C_0 r^{-1+1/N}$.  Let
$1 \leq q \leq 2$, $R \gtrsim r$.
\bi
\item If $\psi$ is a blue wave of arbitrary frequency, we have the finite speed of propagation law
\be{finite-prop} 
\| ((1-P_D) \phi) \psi \|_{L^q(Q(x_D,t_D,C^{-1} r))} \lesssim r^{C-N} E(\phi)^{1/2} E(\psi)^{1/2}
\end{equation}
and the Huygens' principle
\be{phi-eek} \| (P_D \phi) \psi \|_{L^q(Q(x_D,t_D;R) \backslash C^{red}(x_D,t_D;Cr + R^{1/N}))} \lesssim R^{C-N} E(\phi)^{1/2} E(\psi)^{1/2}.
\end{equation}
\item  If $\psi$ is a blue wave of frequency $2^k$ for some $k \geq 0$ such that
$margin(\psi)\gg (2^k r)^{-1+1/N}$, then we have
\be{redblue-h}
 \| (P_D \phi) (P_D \psi) \|_{L^q(Q^{ann}(x_D,t_D;Cr+CR^{1/N},R))} \lesssim R^{C-N} E(\phi)^{1/2} E(\psi)^{1/2}.
\end{equation}
\end{itemize}
\end{lemma}

We will also need to use \eqref{finite-prop}, \eqref{phi-eek}, \eqref{redblue-h} with the roles of $\phi$, $\psi$ and red and blue reversed.  If one repeats the arguments below, one finds that one loses a factor of $2^{Ck}$ on the right-hand sides; however, if $r \geq 2^{k}$ one can absorb these losses into the factors of $r^{C-N}$ or $R^{C-N}$.  Thus we may freely interchange $\phi$ and $\psi$ and red and blue in the $r \geq 2^k$ case. 

\begin{proof}
By H\"older it suffices to verify the claims when $q=2$.

To prove \eqref{finite-prop} we first use \eqref{prop}, \eqref{k-bound}, and \eqref{essential-vanish} to obtain
$$\| (1-P_D) \phi \|_{L^\infty(Q(x_D,t_D;C^{-1} r))} \lesssim r^{C-N} E(\phi)^{1/2},
$$
and \eqref{finite-prop} follows from \eqref{L2-triv}.  A similar argument 
gives
\be{pdp}
 \| P_D \phi \|_{L^\infty(Q(x_D,t_D;R) \backslash C^{red}(x_D,t_D;Cr + R^{1/N}))} \lesssim R^{C-N} E(\phi)^{1/2},
\end{equation}
and \eqref{phi-eek} follows from \eqref{L2-triv}.

We now prove \eqref{redblue-h}.  From \eqref{phi-eek} with $\psi$ replaced by $P_D \psi$, we have
$$ \| (P_D \phi) (P_D \psi) \|_{L^2(Q(x_D,t_D;R) \backslash C^{red}(x_D,t_D;Cr + R^{1/N}))} \lesssim R^{C-N} E(\phi)^{1/2} E(\psi)^{1/2}.$$

By applying \eqref{pdp} with $R$, $r$ replaced by $2^k R$, $2^k r$, and then applying $\T \D_k$, we obtain
$$
 \| P_D \psi \|_{L^\infty(Q(x_D,t_D;R) \backslash C^{blue}(x_D,t_D;Cr + R^{1/N}))} \lesssim R^{C-N} E(\psi)^{1/2},
$$
so by \eqref{L2-triv} we have
$$ \| (P_D \phi) (P_D \psi) \|_{L^2(Q(x_D,t_D;R) \backslash C^{blue}(x_D,t_D;Cr + R^{1/N}))} \lesssim R^{C-N} E(\phi)^{1/2} E(\psi)^{1/2}.
$$
Combining this with the previous estimate we obtain \eqref{redblue-h}.
\end{proof}

The estimate \eqref{redblue-h} is one way of exploiting the transversality of red and blue null directions.  We shall also rely on transversality in other ways, notably in Lemma \ref{bluecone} and in Lemma \ref{mock-local}.

\section{Proof of Proposition \ref{back-prop}}\label{back-sec}

We now have enough machinery to prove Propositions \ref{back-prop} and \ref{r-flip}, except for a somewhat difficult estimate, Lemma \ref{far-cone}, which we defer to Appendix II.

Fix $R \geq 2^{C_0 k}$, and let $\phi$, $\psi$ be red and blue waves of frequency 1 and $2^k$ respectively, obeying the margin requirement \eqref{margin-bound} and the energy normalization \eqref{energy-bound}.  Let $Q_R$ be a cube of side-length $R$. To prove Proposition \ref{back-prop} we have to show that
\be{back-est}
\| \phi \psi \|_{L^p(Q_R)} \leq
(1 - C_0^{-C}) \sup_{2^{C_0 k} \leq \tilde R \leq R; \tilde R^{1/2+4/N} \leq r} A(\tilde R,r,C_0(1+r)) + 2^{CC_0} 2^{\eps k} \K.
\end{equation}
We may of course assume that $\phi$, $\psi$ are nearly extremal in the sense that
\be{extremal}
\| \phi \psi \|_{L^p(Q_R)} \sim A(R).
\end{equation}
We may also assume that 
\be{ar-large}
A(R) \geq 2^{CC_0} \K
\end{equation}
since \eqref{back-est} is trivial otherwise.

Let $0 < \delta < 1/2$ be a small number to be chosen later.  Let $r$ be the supremum of all radii $r \geq C_0 2^{k}$ such that
\be{r-conc}
E_{r, C_0 Q_R}(\phi,\psi) \leq 1-\delta,
\end{equation}
or $r = C_0 2^{k}$ if no such radius exists.

Since $\phi$, $\psi$ are smooth and have finite energy, it is easy to see that $r$ is well defined, and that one can find a disk $D = D(x_0,t_0;r)$ of radius $r$ with $t_0$ in the lifespan of $C_0 Q_R$ for which
\be{energy-c}
\min( \| \phi \|_{L^2(D)}^2, \| \psi \|_{L^2(D)}^2) \geq 1-2\delta.
\end{equation}

Fix the disk $D(x_0,t_0;r)$, and define the disk $D'$ by
$$ D' := C_0^{1/2} D = D(x_0,t_0; C_0^{1/2} r)$$
and define the double conic region $\Omega$ by
$$ \Omega := Q_R \cap C^{purple}(x_0,t_0; C_0(1+r)).$$
We now show that in $Q_R$, the expression $\phi \psi$ is mostly concentrated in $\Omega$.

We split $\phi$ and $\psi$ into red and blue waves respectively as
$$ \phi = (1 - P_{D'})\phi + P_{D'} \phi, \quad \psi = (1 - P_{D'})\psi + P_{D'} \psi $$
where the projection operator $P_{D'}$ is as in the previous section.  

From \eqref{energy-c}, \eqref{energy-bound}, and \eqref{non-local-energy} we have
$$ E((1-P_{D'}) \phi), E((1-P_{D'}) \psi) \lesssim \delta + C_0^{-C}.$$
From \eqref{crude-eq} we thus see that
\be{cooee}
\| ((1-P_{D'}) \phi) (1-P_{D'}) \psi \|_{L^p(Q_R)} \lesssim (\delta + C_0^{-C}) A(R).
\end{equation}
On the other hand, from \eqref{phi-eek} and its analogue with the roles of red and blue reversed (using the assumption $r \geq 2^k$; see the remarks following Lemma \ref{huygens}) we have
$$
\| (P_{D'} \phi) \psi \|_{L^p(Q_R \backslash \Omega)},
\| ((1-P_{D'}) \phi) P_{D'} \psi \|_{L^p(Q_R \backslash \Omega)},
\lesssim C_0^{-C}.$$
Combining these estimates using the triangle inequality we obtain
$$
\| \phi \psi \|_{L^p(Q_R \backslash \Omega)} \lesssim
(\delta + C_0^{-C})A(R).
$$
We thus see from \eqref{extremal} that
$$
\| \phi \psi \|_{L^p(Q_R \backslash \Omega)} \lesssim
(\delta + C_0^{-C}) \| \phi \psi \|_{L^p(Q_R)}.$$
Raising this to the $p^{th}$ power and re-arranging, we see that
\be{midcone}
\| \phi \psi \|_{L^p(\Omega)} \geq 
(1 - C(\delta + C_0^{-C})^p)^{1/p} \| \phi \psi \|_{L^p(Q_R)}.
\end{equation}
We now divide into two cases.  First suppose that we are in the medium or low concentration case $r \geq R^{1/2 + 4/N}$.  In particular we have $r > C_0 2^k$, so from the definition of $r$ we see that \eqref{r-conc} holds.
From the definition of $A(R,r)$ (and the assumption $\delta \ll 1/2$) we thus have
$$ \| \phi \psi \|_{L^p(\Omega)} \leq
(1 - \delta)^{1/p'} A(R,r).$$
Combining these estimates and setting $\delta$ to some small negative power of $C_0$ we obtain the desired estimate \eqref{back-est} as desired (provided that $C_0$ is large enough; \eqref{c0-def} suffices).

Now suppose instead that we are in the highly concentrated case $r \leq C_0 R^{1/2 + 4/N}$.  Define $\tilde R$ by the formula
$$ \tilde R := \max(2^{C_0 k}, (r)^{\frac{1}{1/2 + 4/N}}).$$
thus $2^{C_0 k} \leq \tilde R \leq R$.  In analogy to the previous case, we observe that
$$ \| \phi \psi \|_{L^p(Q(x_0,t_0;\tilde R) \cap \Omega)} \leq
(1 - \delta)^{1/p'} A(\tilde R,r)$$
when $\tilde R > 2^{C_0 k}$.  When $\tilde R = 2^{C_0 k}$ one can use \eqref{ark} (which has already been proven for this value of $\tilde R$) and \eqref{crude-eq} to obtain
$$ \| \phi \psi \|_{L^p(Q(x_0,t_0;\tilde R) \cap \Omega)} \leq 2^{CC_0} 2^{\eps k} \K.$$
In either case, we see that if we could show that
\be{whoops}
\| \phi \psi \|_{L^p(Q_R \backslash Q(x_0,t_0;\tilde R))} \lesssim 
(\delta + C_0^{-C})A(R)
\end{equation}
then we could repeat the argument in the first case and obtain \eqref{back-est} as desired.

It remains to prove \eqref{whoops}.  By \eqref{cooee} and the triangle inequality it suffices to show that
$$
\| (P_{D'} \phi) \psi \|_{L^p(Q_R \backslash Q(x_0,t_0;\tilde R))} \lesssim
(\delta + C_0^{-C})A(R)
$$
and
$$
\| ((1-P_{D'}) \phi) P_{D'} \psi \|_{L^p(Q_R \backslash Q(x_0,t_0;\tilde R))} \lesssim
(\delta + C_0^{-C})A(R).
$$
By a dyadic decomposition and \eqref{ar-large} this will follow from

\begin{lemma}\label{far-cone}
Let $R \geq 2^{C_0 k}$, and let $C_0 2^k < r \leq R^{1/2 + 4/N}$.  Let $D = D(x_0,t_0,C_0^{1/2} r)$ be a disk.  Then for any red wave $\phi$ and blue wave $\psi$ of frequency $1$, $2^k$ respectively with 
$$ \margin(\phi), \margin(\psi) \geq 1/200$$
we have
$$ \| (P_D \phi) \psi \|_{L^p(Q^{ann}(x_0,t_0;R, 2R))},
\| \phi (P_D \psi) \|_{L^p(Q^{ann}(x_0,t_0;R, 2R))}
 \lesssim E(\phi)^{1/2} E(\psi)^{1/2} R^{-1/C}.$$
\end{lemma}

The proof of this lemma shall be postponed to Appendix II, as it requires some machinery which is developed later in this paper, and the arguments used to prove this lemma are not used elsewhere.  The estimate in Lemma \ref{far-cone}
is not of an endpoint type, and could probably be proven by more elementary means than the methods used in the Appendix.

Assuming Lemma \ref{far-cone}, the proof of Proposition \ref{back-prop} is thus complete.
\endprf

\section{Proof of Proposition \ref{r-flip}}\label{flip-sec}

Let $R \geq 2^{C_0 k}$, $0 < c \leq 2^{-C_0}$, and $r > C_0 R$.  Let $\phi$, $\psi$ be red and blue waves of frequency 1 and $2^k$ respectively obeying \eqref{energy-bound} and \eqref{margin-bound}.  Let $Q_R$ be a cube of side-length $R$.  To prove Proposition \ref{r-flip}, it suffices to show that
$$
\| \phi \psi \|_{L^{p}(Q_R)} \leq E_{r,C_0 Q_R}(\phi,\psi)^{1/p}
(1 + Cc) \overline{A(R)}
+ c^{-C} \K.$$
Let $(x_{Q_R}, t_{Q_R})$ be the center of $Q_R$, and let $D$ denote the disk 
$$ D := D(x_{Q_R}, t_{Q_R}; r/2).$$
We then divide
$$ \phi = (1 - P_D)\phi + P_D \phi, \quad \psi = (1 - P_D)\psi + P_D \psi.$$
From Lemma \ref{cutoff} (especially \eqref{local-energy}), Definition \ref{conc-def} and the hypothesis $R \geq 2^{C_0 k}$ we see that $P_D \phi$, $P_D \psi$ obey the relaxed margin requirements \eqref{margin-relaxed} and the energy estimates
$$
E(P_D \phi), E(P_D \psi) \leq E_{r,C_0 Q_R}(\phi,\psi) + C R^{C-N}.$$
From Proposition \ref{upshift} we thus have
$$ 
\| (P_D \phi) (P_D \psi) \|_{L^p(Q_R)} \leq (1+Cc) (E_{r,C_0 Q_R} (\phi,\psi) + C R^{C-N}) \overline{A(R)} + c^{-C} C_0^C \K.$$
The $R^{C-N}$ term can easily be absorbed into the $c^{-C} C_0^C \K$ term by some trivial bound on $\overline{A(R)}$ (Proposition \ref{initial} will do).  By the triangle inequality, we will thus be done if we can show
$$
\| ((1-P_D)\phi) \psi \|_{L^p(Q_R)}, 
\| (P_D\phi) (1-P_D) \psi \|_{L^p(Q_R)}
\leq c^{-C} C_0^C \K.$$
But these estimates follow easily from \eqref{finite-prop} and the counterpart with red and blue reversed (this is legitimate by the hypothesis $R \geq 2^{C_0 k}$ and the remarks following Lemma \ref{huygens}).

\endprf

\section{Energy estimates on light cones of opposite color}\label{opposite-sec}

To prove Theorem \ref{main} it remains only to prove Proposition \ref{lung}.  Of course, we shall not use Propositions \ref{initial}-\ref{r-medium} in the proof of this Proposition.

In this section and the next, we derive two basic estimates needed to prove Proposition \ref{lung}; these are ``pigeonhole-free'' versions of similar estimates in \cite{wolff:cone}.  The first estimate is an extension of \eqref{L2-triv} from a cube to a conic neighbourhood of opposite color.
  The second is a strengthening of \eqref{l2-noloc} when at least one of the waves $\phi$, $\psi$ has low frequency dispersion.  Then, in Section \ref{packet}, we combine these two estimates with a wave packet decomposition to derive a key ingredient (Proposition \ref{fund}) in the proof of Proposition \ref{lung}.

We begin with the energy estimate on cones:

\begin{lemma}\label{bluecone}  Let $\phi$ be a red wave of frequency $2^j$. Then for any $(x_0,t_0) \in \R^{n+1}$ and $R \gtrsim 2^{-j}$ we have
$$ \| \phi \|_{L^2(C^{blue}(x_0,t_0;R))} \lesssim R^{1/2} E(\phi)^{1/2}.$$
\end{lemma}

\begin{proof}
By translation invariance we may take $(x_0,t_0) = 0$, and by scaling we may take $j=0$.

The blue cone is a null surface, and so standard energy estimates will not prove this estimate.  However this can be salvaged since $\phi$ is red, so that the characteristics of $\phi$ will be transverse to the blue cone.

We turn to the details.  Let $U(t)$ be the wave evolution operator defined in \eqref{u-def}.  By \eqref{prop}, it suffices to show that
$$ ( \int \| \chi_{S_t} U(t) f \|_{L^2(\R^n)}^2\ dt)^{1/2} \lesssim R^{1/2} \|f\|_{L^2(\R^n)}$$
for all $f \in L^2(\R^n)$, where $S_t$ is the slice $S_t := \{ x: (x,t) \in C^{blue}(0,0;R)\}$.  

We shall apply the $TT^*$ method.  By duality the above estimate is equivalent to 
$$ \| \int U(t)^* (\chi_{S_t} F(t))\ dt \|_{L^2(\R^n)} \lesssim
R^{1/2} \|F\|_{L^2(\R^{n+1})}.$$
We square this as
$$ \int\int \langle \chi_{S_s} U(s) U(t)^* (\chi_{S_t} F(t)), F(s) \rangle
\ dt ds \lesssim
R \|F\|_2^2$$
where $\langle f,g \rangle = \int_{\R^n} f(x) \overline{g(x)} \ dx$ is the usual
complex inner product.  By Cauchy-Schwarz and Young's inequality it suffices to show that
$$ |\langle \chi_{S_s} U(s) U(t)^* (\chi_{S_t} F(t)), F(s) \rangle|
\lesssim (1 + |s-t|/R)^{-N} \|F(t)\|_2 \|F(s)\|_2$$
for all $s,t$.

When $|s-t| \lesssim R$ this follows from Cauchy-Schwarz and the $L^2$ boundedness of the operator $\chi_{S_s} U(s) U(t)^* \chi_{S_t}$, so we may assume that $|s-t| \gg R$.  
The convolution operator $U(s) U(t)^*$ behaves essentially like $U(s-t)$, and its kernel $K_{s,t}(x)$ satisfies a similar kernel estimate to \eqref{k-bound}, namely
$$ |K_{s,t}(x)| \lesssim (1 + \dist((x,s-t),C^{red}(0,0))^{-N^{10}}.$$
The claim then follows from the transversality of $C^{red}(0,0)$ and $C^{blue}(x_0,t_0)$ and crude estimates.
\end{proof}

Applying $\T$, we see that a similar estimate holds with the roles of red and blue reversed.  As a corollary of both Lemma \ref{bluecone} and its time-reversed counterpart, we obtain the following improvement to \eqref{l1-noloc} when one restricts to a double cone.

\begin{corollary}\label{doublecone}  Let $\phi$, $\psi$ be red and blue waves of frequency $1$, $2^k$ respectively.  Let $R > r \gg 1$, $(x_0,t_0) \in \R^{n+1}$, and $Q_R$ be a cube of side-length $R$.  Then 
$$ \| \phi \psi \|_{L^2(C^{purple}(x_0,t_0;r) \cap Q_R)} \lesssim r^{1/2} R^{1/2} E(\phi)^{1/2} E(\psi)^{1/2}.$$
\end{corollary}

\begin{proof}
The contribution of $C^{blue}(x_0,t_0;r) \cap Q_R$ is acceptable by using Lemma \ref{bluecone} to control $\phi$ and \eqref{L2-triv} to control $\psi$.  Similarly for the contribution of $C^{red}(x_0,t_0;r) \cap Q_R$.
\end{proof}

\section{Bilinear $L^2$ estimates in the low dispersion case.}\label{bilinear-sec}

To prove Proposition \ref{lung} we shall need to study waves with rather low frequency dispersion.  To make this concept precise, we introduce the following notation, which we shall re-use in later sections.

\begin{definition}  If $\phi$ is a red or blue wave (of any frequency), we define the \emph{angular dispersion} of $\phi$ to be the quantity
$$ \diam\{ \frac{\xi}{|\xi|}: (\xi,\tau) \in \supp(\hat \phi) \}.$$
\end{definition}

A key observation (implicit in \cite{mock:cone}, and exploited explicitly in \cite{wolff:cone}) is that there is an improvement to \eqref{l2-noloc} when either $\phi$ or $\psi$ has low dispersion.

\begin{lemma}\label{mock}\cite{mock:cone,wolff:cone} Let $\phi$ be a red wave of frequency 1, and let $\psi$ be a blue wave.  If $\phi$ has angular dispersion $O(1/r)$ for some $r \gg 1$, then
\be{mock-eq} \| \phi \psi \|_{L^2(\R^{n+1})} \lesssim r^{-(n-1)/2} E(\phi)^{1/2} E(\psi)^{1/2}.
\end{equation}
\end{lemma}

\begin{proof}
Let $\psi$ have frequency $2^k$.
From hypothesis, the frequency support of $\phi$ is contained in a
a sector $\Gamma$ of width $O(1/r)$.
By Plancherel's theorem it thus suffices to show that
$$ \| fd\sigma_1 * gd\sigma_2 \|_2 \lesssim r^{-(n-1)/2} \|f\|_{L^2(\Gamma)} \|g\|_{L^2(2^k \Sigma^{blue})}$$
for all $f$ and $g$, where
$d\sigma_1$ and $d\sigma_2$ denote surface measure on $\Gamma$ and $2^k \Sigma^{blue}$ respectively.  By Young's inequality we have
$$ \| fd\sigma_1 * gd\sigma_2 \|_1 \lesssim \|f\|_1 \|g\|_1$$
so it suffices by interpolation to show that
$$ \| fd\sigma_1 * gd\sigma_2 \|_\infty \lesssim r^{-(n-1)} \|f\|_\infty \|g\|_\infty,$$
which by positivity reduces to showing that
\be{sos}
 \| d\sigma_1 * d\sigma_2 \|_\infty \lesssim r^{-(n-1)}.
\end{equation}
But this follows from the observation that the intersection of $\Gamma$
and $x - 2^{k} \Sigma^{blue}$ is always transverse and has $n-1$-dimensional
measure at most $O(r^{-(n-1)})$ for any $x \in \R^{n+1}$.
\end{proof}

This estimate is global in spacetime, however for our purposes it will be convenient to work with a localized form of this estimate. 

\begin{lemma}\label{mock-local} Let $j$ be an integer and $r \gg 1, 2^{-j}$.
Let $\phi$ be a red wave with frequency 1 and angular dispersion $O(1/r)$, and let $\psi$ be a blue wave with frequency $2^j$ such that
$$ \margin(\phi), \margin(\psi) \geq 1/200.$$
Then we have
\be{mock-eq2} \| \phi \psi \|_{L^2(Q)} \lesssim r^{-(n+1)/2} \| \phi \|_{L^2(CQ)} \| \psi \|_{L^2(CQ)} + r^{C-N} E(\phi)^{1/2} E(\psi)^{1/2}
\end{equation}
for all cubes $Q$ of side-length $r$.
\end{lemma}

\begin{proof} 
The idea is to combine Lemma \ref{bluecone} with the localization machinery developed in Section \ref{localize}.

Let $D$ be the disk $D := D(x_Q,r_Q;C^{1/2}r)$.  If $C$ is sufficiently large, then from \eqref{finite-prop} and \eqref{energy-minor} we have
$$ \| ((1 - P_D) \phi) \psi \|_{L^2(Q)} \lesssim r^{C-N} E(\phi)^{1/2} E(\psi)^{1/2}$$
and
$$ \| P_D \phi (1 - P_D) \psi \|_{L^2(Q)} \lesssim r^{C-N} E(\phi)^{1/2} E(\psi)^{1/2}.$$
From the triangle inequality we thus have
$$ \| \phi \psi \|_{L^2(Q)} \leq \| (P_D \phi) (P_D \psi) \|_{L^2(\R^{n+1})} + C r^{C-N} E(\phi)^{1/2} E(\psi)^{1/2}.$$
From Lemma \ref{cutoff} we see that $P_D \phi$ has dispersion $O(1/r)$.
By Lemma \ref{mock} and \eqref{energy-minor} we thus have
$$ \| \phi \psi \|_{L^2(Q)} \lesssim r^{-(n+1)/2} 
\| \phi \|_{L^2(CD)} \| \psi \|_{L^2(CD)}
+ r^{C-N} E(\phi)^{1/2} E(\psi)^{1/2}.$$
The claim then follows by letting $t_D$ range over the lifespan of $Q$, averaging, and applying Cauchy-Schwarz.
\end{proof}

An analogue of this lemma exists with the roles of red and blue reversed, but we shall not invoke this analogue explicitly.

\section{Wave packet decomposition}\label{packet}

Using the results of the previous section, we can now prove the main tool used to derive Proposition \ref{lung}.

\begin{proposition}\label{fund} Let $j,j'$ be integers, and let 
$R \gg 2^{-j}, 2^{j-2j'}$, $0 < c \leq 2^{-C_0}$.  Let $Q$ be a spacetime cube of side-length $R$, $\phi$ be a red wave of frequency $2^j$ with $\margin(\phi) \geq (2^j R)^{-1/2}$, and $\psi$ be a blue wave of frequency $2^{j'}$.

Then there exists a red wave table $\Phi = \Phi_c(\phi,\psi;Q)$ on $Q$ with depth $C_0$, frequency $2^j$ with margin
$$ \margin(\Phi) \geq \margin(\phi) - C (2^j R)^{-1/2}$$ 
such that the following properties hold.
\bi
\item ($[\Phi]_{C_0}$ approximates $\phi$) We have 
\be{non-conc2}
\| (|\phi| - [\Phi]_{C_0}) \psi \|_{L^2(I^{c,C_0}(Q))} \lesssim c^{-C} 
(\frac{2^j}{R})^{(n-1)/4} E(\phi)^{1/2} E(\psi)^{1/2}.
\end{equation}
\item (Bessel inequality) We have
\be{bessel}
E(\Phi) \leq (1+ Cc) E(\phi)
\end{equation}
\item (Persistence of non-concentration)  For any $r \gtrsim R^{1/2 + 1/N} 2^{-j(1/2 - 1/N)}$, we have
\be{conc-persist}
E_{r(1 - C_0 (2^j r)^{-1/2N}),C_0 Q}(\Phi, \psi) \leq (1 + Cc) E_{r,C_0 Q}(\phi,\psi) 
\end{equation}
\end{itemize}
\end{proposition}

Roughly speaking, $\Phi^{(q)}$ is the portion of $\phi$ which concentrates in $q$.  The Proposition is phrased using frequencies $2^j$ and $2^{j'}$ rather than 1 and $2^k$ in order to make it symmetric in $\phi$ and $\psi$, as we shall also need the time reversal of Proposition \ref{fund}.

\begin{proof}  The quantities $r,R$ have the units of length, while the frequencies $2^j$, $2^{j'}$ have the units of inverse length.  Finally, $c$ and the margins are dimensionless.  One can then verify that the entire Proposition is dimensionally correct and thus scale invariant.  By scaling we may thus set $j=0$.  By translation invariance we may assume $Q$ is centered at the origin.  Set $r := 2^{-J} R$, where $J$ is chosen so that $r \sim \sqrt{R}$.

We may assume that $c > R^{-1/10n}$ (for instance), since otherwise the claim is trivial by setting $\Phi = 0$ (and using some trivial bound such as \eqref{l2-noloc} to show \eqref{non-conc2}.

Let $\Ecal$ be a maximal $1/r$-separated subset of
$S^{n-1} \cap \underline \Sigma$, and let $L$ denote the lattice $L := c^{-2} r \Z^n$.
We define a \emph{red tube} to be any set $T = T(\omega_T,x_T)$ of the form
$$ T = \{ (x,t) \in \R^{n+1}: |x - (x_T + \omega_T t)| \leq r \}$$
where $x_T \in L$ and $\omega_T \in \Ecal$.  We let $\Tcal = \Tcal^{red}$ denote the set of all red tubes.
If $T$ is a red tube, we define the cutoff function $\tilde \chi_T$ on $\R^{n+1}$ by
\be{tchit-def} \tilde \chi_T(x,t) := \tilde \chi_{D(x_T + \omega_T t,t;r)}(x). \end{equation}

We shall need the following careful decomposition of $\phi$ into wave packets which are concentrated on tubes in $\Tcal$.

\begin{lemma}\label{tube-decomp}  With the above notation, one can find for each $T \in \Tcal$ a red wave $\phi_T$ with frequency 1 and angular dispersion at most $C R^{-1/2}$, such that
\bi
\item We have the margin estimate
\be{margin-phit}
\margin(\phi_T) \geq \margin(\phi) - C R^{-1/2}.
\end{equation}
\item The map $\phi \mapsto \phi_T$ is linear for each $T$, and 
\be{recon-1}\phi = \sum_{T \in \Tcal} \phi_T.
\end{equation}
\item We have
\be{local2} E(\phi_T) \lesssim c^{-C} \| \tilde \chi_T(t) \phi(t) \|_2
\end{equation}
for all $T \in \Tcal$ and $t$ in the lifespan of $C_0 Q$.  
\item If $\dist(T,Q) \geq C_0 R$ then
\be{local3} \| \phi_T \|_{L^\infty(Q)} \lesssim \dist(T,Q)^{C-N} E(\phi)^{1/2}.
\end{equation}
\item We have
\be{local} 
\sum_{T \in \Tcal} \sup_{q \in \Qcal_J(Q)} 
\tilde \chi_T(x_q,t_q)^{-3}
 \| \phi_T \|_{L^2(Cq)}^2
\lesssim c^{-C} r E(\phi).
\end{equation}
\item We have the Bessel inequality
\be{bessel-T} (\sum_{q_0} E(\sum_{T \in \Tcal} \m_{q_0,T} \phi_T))^{1/2} \leq (1 + Cc) E(\phi)^{1/2}
\end{equation}
whenever $q_0$ runs over a finite index set and the $\m_{q_0,T}$ are non-negative numbers such that
\be{eta-sum}
\sum_{q_0} \m_{q_0,T} = 1
\end{equation}
 for all $T \in \Tcal$.
\end{itemize}
\end{lemma}

Roughly speaking, $\phi_T$ is the portion of $\phi$ which has frequency support in the sector of width $1/r$ and direction $(\omega_T,1)$, and is spatially concentrated in $T$.  A naive microlocalization to this region of space and frequency, taking some care to ensure that the $\phi_T$ are still waves, would obtain most of the above properties, but would probably need to replace the $(1+Cc)$ factor in \eqref{bessel-T} by a larger constant, which would then cause a similar unacceptable loss in \eqref{energy-est} and then destroy the induction.  This necessitates a delicate construction of the $\phi_T$ based on averaging.  Because the details of the proof are technical and not particularly relevant to the rest of the argument, we defer the proof of this Lemma to Appendix I, and continue with the proof of Proposition \ref{fund}.

Using Lemma \ref{tube-decomp}, we can now define $\Phi$ by

\be{phi-Om-def} \Phi^{(q_0)} := \sum_{T \in \Tcal} \frac{\m_{q_0,T}}{\m_T} \phi_T
\end{equation}
for all $q_0 \in \Qcal$, where
\be{eta-def} \m_{q_0,T} := \| \psi \tilde \chi_T \|_{L^2(q_0)}^2 
+ R^{-10n} E(\psi),
\end{equation}
and
\be{etat-def} \m_T := \sum_{q_0 \in \Qcal_{C_0}(Q)} \m_{q_0,T}
= \| \psi \tilde \chi_T \|_{L^2(Q)}^2 + R^{-10n} 2^{(n+1)C_0} E(\psi).
\end{equation}
The $R^{-10n} E(\psi)$ factor is only present in \eqref{eta-def} to ensure that $\m_T$ does not completely degenerate to zero.  One can think of $\Phi^{(q_0)}$ as consisting of those wave packets $\phi_T$ such that $\psi \chi_T$ concentrates in $q_0$.

It is clear that $\Phi$ is a red wave with 
$$ \margin(\Phi) \geq \margin(\phi) - CR^{-1/2}$$
and that 
\be{recon}
 \phi = \sum_{q_0 \in \Qcal_{C_0}(Q)} \Phi^{(q_0)}.
\end{equation}
The estimate \eqref{bessel} follows from \eqref{bessel-T}.

We now show \eqref{conc-persist}.  Let $D = D(x_0,t_0; r(1 - C_0 r^{-1/2N}))$ be any disk of radius $r(1 - C_0 r^{-1/2N})$ with $t_0$ in the lifespan of $C_0Q$.  By \eqref{bessel} it then suffices to show that
\be{chico}
\| \Phi \|_{L^2(D)} \| \psi\|_{L^2(D)}
\leq (1 + Cc) E_{r',C_0 Q}(\phi,\psi).
\end{equation}
Let $D'$ denote the slightly larger disk $D' := D(x_0,t_0,r - \frac{C_0}{2} r^{-1/2N})$, and let $D''$ denote the even larger disk $D'' := D(x_0, t_0, r)$.  We may divide 
$$ \phi = P_{D'} \phi + (1-P_{D'}) \phi =: \phi_1 + \phi_2.$$
The map $\phi \mapsto \Phi$ is linear, so we may write $\Phi =: \Phi_1 + \Phi_2$ accordingly.  From \eqref{bessel} and \eqref{local-energy} we have
\be{phi1-bound}
\| \Phi_1 \|_{L^2(D)} \leq E(\Phi_1)^{1/2} \leq (1+Cc) E(\phi_1)^{1/2}
\leq (1 + Cc) \| \phi \|_{L^2(D'')} + r^{C-N} E(\phi).
\end{equation}
Next, we claim that
\be{phi2-bound}
\| \Phi_2 \|_{L^2(D)} \lesssim r^{C-N} E(\phi).
\end{equation}
To see this,  we first consider the tubes $T$ which do not intersect $CQ$.  By \eqref{phi-Om-def}, \eqref{local3}, and H\"older, their contribution is acceptable.  Thus we need only consider those tubes which intersect $CQ$.  By the triangle inequality it suffices to show 
$$ \| (\phi_2)_T \|_{L^2(D)} \lesssim r^{C-N} E(\phi)$$
for each tube $T$.  

First suppose that $\dist(T,D) \leq R^{1/2 + 1/100N}$.  Then by  \eqref{essential-vanish}, \eqref{energy-minor}, and the assumption $r \gtrsim R^{1/2 + 1/N}$ we see that
$$ \|\tilde \chi_T(t_0) \phi_2(t_0) \|_2 \lesssim r^{C-N} E(\phi).$$
The claim then follows from \eqref{local2}.

Now suppose that $\dist(T,D) \geq R^{1/2 + 1/100N}$.  By \eqref{local}, the hypothesis $c > R^{-1/10}$ and Bernstein's inequality (or Sobolev embedding) we see that
$$ \| (\phi_2)_T \|_{L^\infty(q)} \lesssim r^{C-N} E(\phi_2) \lesssim r^{C-N} E(\phi)$$
whenever $\dist(q,T) \geq R^{1/2 + 1/200N}$.  The claim then follows.  This completes the proof of \eqref{phi2-bound}.

Combining \eqref{phi1-bound} and \eqref{phi2-bound}, we see that
$$ \| \Phi \|_{L^2(D)} \leq (1 + Cc) \| \phi \|_{L^2(D'')} + r^{C-N} E(\phi),$$
and \eqref{chico} easily follows (using the trivial inequality $\| \psi \|_{L^2(D)} \leq \| \psi \|_{L^2(D'')}$).  This proves \eqref{conc-persist}.

We now turn to the proof of \eqref{non-conc2}. By \eqref{recon} we have
$$ (|\phi| - [\Phi]_{C_0}) |\psi| \leq \sum_{q_0 \in \Qcal_{C_0}(Q)} |\Phi^{(q_0)} \psi| (1 - \chi_{q_0}),$$
so by the triangle inequality it suffices to show that
\be{conc-p}
\| \Phi^{(q_0)} \psi \|_{L^2(I^{c,C_0}(Q) \backslash q_0)} \lesssim c^{-C}
r^{-(n-1)/2} E(\phi)^{1/2} E(\psi)^{1/2}.
\end{equation}
for each $q$ in $\Qcal_{C_0}(Q)$.

Fix $q_0$.  
We shall use a (heavily disguised) version of the arguments in Wolff
\cite{wolff:cone}.  If $q \in \Qcal_J(Q)$ intersects $I^{c,C_0}(Q) \backslash q_0$, then $\dist(q,q_0) \gtrsim cR$.  By squaring \eqref{conc-p}, we thus reduce to showing that
\be{red}
\sum_{q \in \Qcal_J(Q): \dist(q,q_0) \gtrsim cR}
\| \Phi^{(q_0)} \psi \|_{L^2(q)}^2 
\lesssim c^{-C} r^{-(n-1)} E(\phi) E(\psi).
\end{equation}
Consider a single summand from \eqref{red}.
From \eqref{phi-Om-def} and the triangle inequality
we have
\be{turn}
\| \Phi^{(q_0)} \psi \|_{L^2(q)}
\leq \sum_{T \in \Tcal} \frac{\m_{q_0,T}}{\m_T}
\| \phi_T \psi \|_{L^2(q)}.
\end{equation}
Consider the tubes $T$ which do not intersect $CQ$.  By \eqref{local3} and \eqref{L2-triv}, their total contribution to \eqref{turn} is $O(R^{C-N} E(\phi)^{1/2} E(\psi)^{1/2})$.  As for the tubes $T$ which do intersect $CQ$,
we may apply Lemma \ref{mock-local} (since $r \gg 1, 2^{-j'}$) to obtain
$$ \| \phi_T \psi \|_{L^2(q)}
\lesssim r^{-(n+1)/2} 
\| \psi \|_{L^2(Cq)}
\| \phi_T \|_{L^2(Cq)}
+ R^{C-N} E(\phi_T)^{1/2} E(\psi)^{1/2}.$$
The total contribution of the error term to \eqref{turn} is $O(R^{C-N} E(\phi)^{1/2} E(\psi)^{1/2})$ by \eqref{bessel-T} and the fact that there are only $O(R^C)$ tubes being summed here.  Combining all these estimates we obtain
$$
\| \Phi^{(q_0)} \psi \|_{L^2(q)}
\lesssim 
r^{-(n+1)/2} \| \psi \|_{L^2(Cq)}
\sum_{T \in \Tcal} \frac{\m_{q_0,T}}{\m_T} \| \phi_T \|_{L^2(Cq)}
+ R^{C-N} E(\phi)^{1/2} E(\psi)^{1/2}.$$
Inserting this back into \eqref{red}, we see that it suffices to show that
\be{new-targ}
\sum_{q \in \Qcal_J(Q): \dist(q,q_0) \gtrsim cR}
\| \psi \|_{L^2(Cq)}^2
(\sum_{T \in \Tcal} \frac{\m_{q_0,T}}{\m_T} 
\| \phi_T \|_{L^2(Cq)})^2
\lesssim c^{-C} r^2 E(\phi) E(\psi).
\end{equation}
Using the trivial estimate
$$ \frac{\m_{q_0,T}}{\m_T} \leq \frac{\m_{q_0,T}^{1/2}}{\m_T^{1/2}}$$
followed by Cauchy-Schwarz, we have
$$
(\sum_{T \in \Tcal} \frac{\m_{q_0,T}}{\m_T} \| \phi_T \|_{L^2(Cq)})^2
\leq
(\sum_{T \in \Tcal} \frac{ \| \phi_T\|_{L^2(Cq)}^2 }{\m_T \tilde \chi_T(x_q,t_q)})
(\sum_{T \in \Tcal} \m_{q_0,T} \tilde \chi_T(x_q,t_q)).
$$
By \eqref{eta-def}, we have
$$
\sum_{T \in \Tcal} \m_{q_0,T} \tilde \chi_T(x_q,t_q) \lesssim \| \psi \chi\|_2^2 + R^{-5n} E(\psi)
$$
where
$$ \chi := (\sum_{T \in \Tcal} \tilde \chi_T(x_q,t_q) \tilde \chi_T^2)^{1/2}
\chi_{q_0}.$$
Since $\dist(q,q_0) \gtrsim cR$, we see from elementary geometry that
$$ \chi(x,t) \lesssim c^{-C} (1 + \frac{\dist((x,t),C^{red}(x_q,t_q))} {r})^{-10n}.$$
From Lemma \ref{bluecone} and the triangle inequality we thus have
$$ \| \psi \chi\|_2^2 \lesssim c^{-C} r E(\psi).$$
Inserting all these estimates back into \eqref{new-targ}, we see that it will suffice to show that
\be{job}
\sum_{q \in \Qcal_J(Q)}
\| \psi \|_{L^2(Cq)}^2
\sum_{T \in \Tcal} \frac{ \| \phi_T \|_{L^2(Cq)}^2 }
{ \m_T \tilde \chi_T(x_q,t_q) }
\lesssim c^{-C} r E(\phi).
\end{equation}
We re-arrange the left-hand side as
$$ \sum_{T \in \Tcal} \sum_{q \in \Qcal_J(Q)} 
\tilde \chi_T(x_q,t_q)^{-3} \| \phi_T \|_{L^2(Cq)}^2 
\frac{\| \psi \|_{L^2(Cq)}^2 \tilde \chi_T(x_q,t_q)^2}
{ \m_T }
$$
and estimate this by
$$ \sum_{T \in \Tcal} (\sup_{q \in \Qcal_J(Q)} 
\tilde \chi_T(x_q,t_q)^{-3} \| \phi_T \|_{L^2(Cq)}^2)
\sum_{q \in \Qcal_J(Q)} \frac{\| \psi \|_{L^2(Cq)}^2 \tilde \chi_T(x_q,t_q)^2}
{ \m_T }
$$
Since the inner sum is $O(1)$ by \eqref{etat-def}, the desired estimate \eqref{job} then follows from \eqref{local}.  This proves \eqref{non-conc2}.
\end{proof}

We define the wave table $\Psi_{c}(\phi,\psi;Q)$ in analogy to $\Phi_{c}(\phi,\psi;Q)$ by time reversal:
$$ \Psi_{c}(\phi,\psi;Q) := \T \Phi_{c}(\T\psi, \T\phi;\T Q).$$
Of course, one has the analogue of Proposition \ref{fund} for $\Psi$ but with the roles of red and blue reversed.

\section{Proof of Proposition \ref{lung}}\label{fundamental}

We are now ready to prove Proposition \ref{lung}, in which $\phi \psi$ is replaced by a quilted analogue $[\Phi]_k [\Psi]_{C_0}$.  

We begin with the construction of $\Phi$; this will be achieved by iterating
Proposition \ref{fund}.

We define recursively the wave tables $\phi_j$ on $Q$ of depth $j$ for all $0 \leq j \leq k$ which are multiples of $C_0$ (recall that we assumed $k$ to be a multiple of $C_0$ in Section \ref{top-sec}).  Set $\phi_0$ to be the wave table on $Q$ of depth 0
$$ (\phi_0)^{(Q)} := \phi$$ 
and then define inductively
$$ \phi_{j+C_0}^{(q)} := \Phi_{c 2^{-(k-j)/N}}(\phi_j^{(q)}, \psi; q)$$
for all $q \in \Qcal_j(Q)$.  We then choose $\Phi$ to be $\Phi := \phi_k$.

By induction we see that 
$$ \margin(\phi_j) \geq 1/100 - (2^k/R)^{1/N} - C (2^j/R)^{1/2}$$ 
for each $j$.  This gives the desired margin requirements on $\Phi$ since
$R \gg 2^k$.  From \eqref{bessel} we have
$$
E(\phi_j) \leq (1 + Cc 2^{-(k-j)/N} ) E(\phi_{j-1})
$$
for all ${C_0} < j \leq k$.  Telescoping this, we obtain
\be{bessel-j}
E(\phi_j) \leq (1 + Cc) E(\phi) = 1 + Cc
\end{equation}
for all $j$.
Thus \eqref{energy-est} holds for $\Phi$.  A similar iteration of \eqref{conc-persist} (absorbing any factor of $k$ which appears into a small power of $r^{1/N})$ yields
\be{step-1}
E_{r(1 - r^{-1/3N}),C_0 Q}(\Phi, \psi) \leq (1 + Cc) E_{r,C_0 Q}(\phi,\psi)
\end{equation}
for all $r \gtrsim R^{1/2 + 2/N}$.

Let $C_0 \leq j < k$ be a multiple of $C_0$, and let $q^*$ be a cube in $\Qcal_j(Q)$.
From \eqref{non-conc2} with $R$ replaced by $2^{1-j}R$, 
$Q$ replaced by $q^*$, and $\phi$ replaced by $\phi_j^{(q^*)}$, we see that
$$
\| ([\phi_j]_j - [\phi_{j+C_0}]_{j+C_0}) \psi \|_{L^2(X(Q) \cap q^*)}
\lesssim
c^{-C} 2^{C(k-j)/N}
(2^{1-j} R)^{-\frac{n-1}{4}} E(\phi_j^{(q^*)})^{1/2}.$$
Square-summing this in $q^*$ and using \eqref{bessel-j} we obtain
\be{phijj}
\| ([\phi_j]_j - [\phi_{j+C_0}]_{j+C_0}) \psi \|_{L^2(X(Q))}
\lesssim
c^{-C} 2^{C(k-j)/N} (2^{-j} R)^{-\frac{n-1}{4}}.
\end{equation}
On the other hand, from Lemma \ref{phipsi} with $R$ replaced by $2^k R$ and $\phi$ replaced by $\phi_j$, we have
$$ \| [\phi_j]_j \psi \|_{L^1(X(Q))}
 \lesssim 2^{-j/2} R.$$
Replacing $j$ by $j+C_0$ and then subtracting, we obtain
$$ \| ([\phi_j]_j - [\phi_{j+C_0}]_{j+C_0}) \psi \|_{L^1(X(Q))}
 \lesssim 2^{CC_0} 2^{-j/2} R.$$
Interpolating this with \eqref{phijj} using \eqref{p-def}, we obtain
$$ \| ([\phi_j]_j - [\phi_{j+C_0}]_{j+C_0}) |\psi| \|_{L^p(X(Q))}
 \lesssim c^{-C} 2^{C(k-j)/N} 2^{-j(\frac{1}{p}-\frac{1}{2})}
= c^{-C} 2^{-\alpha(k-j)} \K$$
for some constant $\alpha > 0$.  Telescoping this with the triangle inequality, we obtain
\be{dyad}
\| \phi \psi \|_{L^p(X(Q))} \leq \| [\Phi]_k \psi \|_{L^p(X(Q))} 
+ c^{-C} \K.
\end{equation}

Having constructed $\Phi$, we now define $\Psi$ as
$$ \Psi := \Psi_c(\Phi,\psi;Q).$$

The estimate \eqref{energy-est} for $\Psi$ follows from the time reversal of \eqref{bessel}. It remains to prove \eqref{triad}, \eqref{triad-cone}, and \eqref{energy-persist}.

To prove \eqref{triad}, we observe 
from the analogue of \eqref{non-conc2} that
$$
\| \Phi (|\psi| - [\Psi]_{C_0}) \|_{L^2(X(Q))} \lesssim c^{-C} 2^{C{C_0}} 
(R/2^k)^{-(n-1)/4} E(\Phi)^{1/2} E(\psi)^{1/2}.
$$
From \eqref{energy-est} and \eqref{quilt-est} we thus have
\be{tri-2}
\| [\Phi]_k (|\psi| - [\Psi]_{C_0}) \|_{L^2(X(Q))} \lesssim c^{-C} 2^{C{C_0}} 
(R/2^k)^{-(n-1)/4}.
\end{equation}
From Lemma \ref{phipsi} we have
\be{phip}
\| [\Phi]_k \|_{L^2(X(Q))} \lesssim 2^{-k} R^{1/2}.
\end{equation}
From \eqref{energy-est} we thus we have
$$\| [\Phi]_k \psi \|_{L^1(X(Q))}, \| [\Phi]_k \Psi \|_{L^1(X(Q))} \lesssim 2^{-k/2} R
$$
so from \eqref{quilt-est} and the triangle inequality we have
$$ 
\| [\Phi]_k (|\psi| - [\Psi]_{C_0}) \|_{L^1(X(Q))} \lesssim 2^{-k/2} R.$$
Interpolating this with \eqref{tri-2}, we obtain
\be{dyad-2}
\| [\Phi]_k (|\psi| - [\Psi]_{C_0}) \|_{L^p(X(Q))} \lesssim c^{-C} 2^{C{C_0}} \K,
\end{equation}
and \eqref{triad} follows from this, \eqref{dyad}, and the triangle inequality.

Now we prove \eqref{triad-cone}. We shall assume that $2^k r < R$ since \eqref{triad-cone} follows from \eqref{triad} otherwise.  For brevity let $\Omega$ denote the region
$$ \Omega := X(Q) \cap C^{purple}(x_0,t_0;r).$$
We will need to localize \eqref{dyad} and \eqref{dyad-2} to $\Omega$.  We begin with the localization of \eqref{dyad}.  By \eqref{quilt-est} and Corollary \ref{doublecone} we have
$$ \| [\phi_j]_j \psi \|_{L^1(\Omega)}
\leq
\| \phi_j \psi \|_{L^1(\Omega)}
 \lesssim r^{1/2} R^{1/2} \leq (2^k r/R)^{1/2} 2^{-j} R.$$
Replacing $j$ by $j+C_0$ and then subtracting, we obtain
$$ \| ([\phi_j]_j - [\phi_{j+C_0}]_{j+C_0}) \psi \|_{L^1(\Omega)}
 \lesssim (2^k r/R)^{1/2} 2^{-j} R.$$
Interpolating this with \eqref{phijj} using \eqref{p-def}, we obtain
$$ \| ([\phi_j]_j - [\phi_{j+C_0}]_{j+C_0}) \psi \|_{L^p(\Omega)}
 \lesssim c^{-C} 2^{C(k-j)/N} 2^{-j(\frac{1}{p}-\frac{1}{2})}
(2^k r/R)^{1/N}.$$
Telescoping this with the triangle inequality, we obtain
\be{dyad-alt}
\| \phi \psi \|_{L^p(\Omega)} \leq \| [\Phi]_k \psi \|_{L^p(\Omega)} 
+ C c^{-C} 2^{C{C_0}} \K (2^k r/R)^{1/N}.
\end{equation}
Now we localize \eqref{dyad-2}.  From \eqref{quilt-est} and Corollary \ref{doublecone} we have
$$\| [\Phi]_k \psi \|_{L^1(\Omega)}, \| [\Phi]_k \Psi \|_{L^1(\Omega)} \lesssim r^{1/2} R^{1/2} = (2^k r/R)^{1/2} 2^{-k/2} R.$$
so from \eqref{quilt-est} and the triangle inequality we have
$$ 
\| [\Phi]_k (|\psi| - [\Psi]_{C_0}) \|_{L^1(\Omega)} \lesssim (2^k r/R)^{1/2} 2^{-k/2} R.$$
Interpolating this with \eqref{tri-2}, we obtain
$$
\| [\Phi]_k (|\psi| - [\Psi]_{C_0}) \|_{L^p(\Omega)} \lesssim (2^k r/R)^{1/N} c^{-C} 2^{C{C_0}} \K,
$$
and \eqref{triad} follows from this, \eqref{dyad}, and the triangle inequality.

Finally, we show \eqref{energy-persist}.  From the time-reversed version of
\eqref{conc-persist}, and conceding various powers of $2^k$ and $r^{1/N}$, we see that
$$
E_{r(1 - r^{-1/3N}),C_0 Q}(\Phi, \Psi) \leq (1 + Cc) E_{r(1 - r^{1/2N}),C_0 Q}(\Phi,\psi)
$$
for all $r \gtrsim R^{1/2 + 3/N}$.  Combining this with \eqref{step-1} we obtain \eqref{energy-persist} as desired.  This completes the proof of Proposition \ref{lung}, and hence of Theorem \ref{main}.
\endprf

\section{Null form estimates}\label{final}

In this section we now apply Theorem \ref{main} to obtain some nearly sharp null form estimates.  The arguments in this section do not require the methods used above to prove Theorem \ref{main}.

Let $D_0$, $D_+$, $D_-$ denote the Fourier multipliers
\bas
\widehat{D_0 \phi}(\xi,\tau) &:= |\xi| \hat \phi(\xi,\tau)\\
\widehat{D_+ \phi}(\xi,\tau) &:= (|\xi| + |\tau|) \hat \phi(\xi,\tau)\\
\widehat{D_- \phi}(\xi,\tau) &:= \left| |\xi| - |\tau| \right| \hat \phi(\xi,\tau).
\end{align*}

If $\phi$ is a solution to the free wave equation and $s \in \R$, we define $\phi[0]$ to be the vector $(\phi(0), D_0^{-1} \phi_t(0))$.  In particular, we have
$$ \| \phi[0] \|_{\dot H^s} = ( \| \phi(0) \|_{\dot H^s}^2 + \| \phi_t(0) \|_{\dot H^{s-1}}^2)^{1/2}$$
where $\dot H^s := D_0^{-s} L^2$ is the homogeneous Sobolev space of order $s$.

In \cite{damiano:null} the following problem was considered:

\begin{problem}\label{problem} Determine the set of all exponents $(p, \beta_0, \beta_+, \beta_-, \alpha_1, \alpha_2)$ such that one has the estimates
\be{null} \| D_0^{\beta_0} D_+^{\beta_+} D_-^{\beta_-} (\phi \psi) \|_{L^p(\R^{n+1})} \lesssim \| \phi[0] \|_{\dot H^{\alpha_1}} \| \psi[0] \|_{\dot H^{\alpha_2}}
\end{equation}
for all vector-valued solutions $\phi$, $\psi$ to the free wave equation (not necessarily red or blue).
\end{problem}

This problem was resolved in \cite{damiano:null} for $p=2$, but is largely open otherwise (with some partial results in \cite{tv:cone2}, \cite{kl-tar:yang-mills}).  A successful resolution to this problem (and of its generalization to mixed Lebesgue norms $L^q_t L^r_x$ and to when $\phi$, $\psi$ solve the inhomogeneous wave equation) is likely to have application to the low-regularity behaviour of non-linear wave equations.

To see the connection between Problem \ref{problem} and Theorem \ref{main}, we first observe the consequence of Theorem \ref{main}, which addresses Problem \ref{problem} in the frequency-localized case.

\begin{proposition}\label{toy}  Let $\beta$ be a real number, and $k,l \geq 0$ be integers.  Let $\phi$ and $\psi$ be waves which have frequency supports in the sectors
\be{phi-sec} \{ (\xi_1, \xi_2, \xi'', \tau): \tau = |\xi|, \xi_1 \sim 1, \xi_2 \sim 2^{-l}, |\xi''| \lesssim 2^{-l} \}
\end{equation}
and
\be{psi-sec}
 2^k \{ (\xi_1, \xi_2, \xi'', \tau): \tau = -|\xi|, -\xi_1 \sim 1, \xi_2 \sim 2^{-l}, |\xi''| \lesssim 2^{-l} \}
\end{equation}
respectively.  Then we have
\be{est-d} \| |\Box|^{\beta} (\phi \psi) \|_p, \| |\Box|^\beta (\phi \overline \psi) \|_p \lesssim 
2^{\beta (k-2l)} 2^{l(\frac{n+1}{p} - (n-1))}
2^{k(\frac{1}{p}-\frac{1}{2}+\eps)} E(\phi)^{1/2} E(\psi)^{1/2}
\end{equation}
for all $2 \geq p \geq p_0$, $\eps > 0$, where $|\Box|$ is the multiplier $|\Box| := D_+ D_-$, the energy $E(\cdot)$ is defined as in \eqref{energy-def}, and the implicit constants may depend on $\beta$, $\eps$.
\end{proposition}

\begin{proof}
Consider the Minkowski-conformal linear transformation $L: \R^{n+1} \to \R^{n+1}$ given by
$$ L(\xi_1,\xi',\tau) := (\frac{\xi_1 + \tau}{2} + 2^{-2l} \frac{\xi_1 - \tau}{2}, 2^{-l}\xi', \frac{\xi_1 + \tau}{2} - 2^{-2l} \frac{\xi_1 - \tau}{2}),$$
and define the associated operator $T_L$ by
$$ T_L \phi := \phi \circ L^*,$$
where $L^*$ is the adjoint of $L$.
A routine computation shows that $T_L |\Box| = 2^{-2j} |\Box| T_L$, so that
$$ T_L(|\Box|^\beta \phi \psi) = 2^{-2j\beta} |\Box|^\beta (T_L \phi) 
(T_L \psi).$$
Also, we observe from Plancherel that $E(T_L \phi) \sim 2^{(n-1)j} E(\phi)$ and $E(T_L \psi) \sim 2^{(n-1)j} E(\psi)$.  Finally, we have $\| T_L F \|_p =
2^{-(n+1)l/p} \|F\|_p$.  Combining all these facts we see that to prove \eqref{est-d} it suffices to do so when $l=0$.

Set $l=0$.  We can write
$$ |\Box|^{\beta}(\phi \psi)(x,t) = 2^{nk} 2^{\beta k} \int\int e^{2\pi i x\cdot (2^k \eta + \xi)} e^{2\pi i t (|\xi| - 2^k |\eta|)}
m(\xi,\eta)
\hat f(\xi) \hat g(2^k \eta)\ d\xi d\eta$$
where $f(x) := \phi(x,0)$ and $g(x) := \psi(x,0)$, and $m$ is the symbol
$$ m(\xi,\eta) := C a(\xi) a(\eta)  
(|\eta| |\xi| + \langle \eta,\xi \rangle)^{\beta}.$$

Since $m$ is smooth and compactly supported, we may decompose $m$ as a Fourier series
$$ m(\xi,\eta) = \sum_{j,j' \in L} c_{j,j'} e^{-2 \pi i j \cdot \xi} e^{-2\pi i j' \cdot \eta}$$
for $|\xi|, |\eta| \lesssim 1$, where $L$ is some discrete lattice and $c_{j,j'}$ are rapidly decreasing co-efficients (uniformly in $k$).  This implies that
$$ D_-^{\beta_-}(\phi \psi) = \sum_{j,j' \in L} c_{j,j'} \phi_j \psi_{j'}$$
where
$$ \phi_j(x,t) := \phi(x-j,t), \quad \psi_{j'}(x,t) := \psi(x - 2^{-k} j',t).$$
The first claim of \eqref{est-d} then follows from this decomposition, the triangle inequality, \eqref{est}, and the observation that $\phi_j$, $\psi_{j'}$ have the same energy as $\phi$, $\psi$ respectively.  The second claim is proven similarly but also uses the observation that $\phi_j \overline{\psi_{j'}}$ has the same magnitude as $\phi_j \psi_{j'}$.
\end{proof}

We make the technical remark that the above estimates continue to hold if one replaces $\beta$ by the complex number $\beta + it$, with constants that grow at most exponentially in $t$.

We now consider the general setting of Problem \ref{problem}, where no frequency restrictions are assumed on $\phi$ or $\psi$.

By considering a large number of key examples, it was shown in \cite{damiano:null} that the conditions
\ba
\beta_0 + \beta_+ + \beta_- &= \alpha_1 + \alpha_2 + \frac{n+1}{p} - n\label{scaling}\\
p &\geq p_0\label{p-cond}\\
\beta_- &\geq \frac{n+1}{2p} - \frac{n-1}{2} \label{beta-minus}\\
\beta_0 &\geq \frac{n+1}{p} - n \nonumber\\
\beta_0 &\geq \frac{n+3}{p} - (n+1)\label{beta-0}\\
\alpha_1 + \alpha_2 &\geq \frac{1}{p} \label{a12-1}\\
\alpha_1 + \alpha_2 &\geq \frac{n+3}{p} - n \label{a12-2}\\
\alpha_i &\leq \beta_- + \frac{n}{2}\nonumber\\
\alpha_i &\leq \beta_- + \frac{n-1}{2} + \frac{n+1}{2}(\frac{1}{2}-\frac{1}{p})\nonumber\\
\alpha_i &\leq \beta_- + \frac{n-1}{2} + (n+2)(\frac{1}{2}-\frac{1}{p})\label{alpha-i}
\end{align}
were necessary for \eqref{null}, where $i=1,2$.  In that paper it was also shown that the above conditions were also necessary when $p=2$, except for the endpoint cases when one has equality in \eqref{beta-0}, [\eqref{beta-minus} and \eqref{a12-1}], or [\eqref{beta-minus} and \eqref{alpha-i} for some $i=1,2$], with \eqref{null} being false in these cases.  It was then conjectured that for the above conditions are similarly necessary when $p\neq 2$, except perhaps for some endpoints.

By combining Proposition \ref{toy} with dyadic decomposition arguments as in
\cite{tv:cone2}, we can obtain the following progress on this conjecture in the $p \leq 2$ case.

\begin{theorem}\label{n}  Let $p_0 \leq p \leq 2$.  Then 
\eqref{null} holds whenever \eqref{scaling} holds, \eqref{beta-minus}, \eqref{beta-0}, \eqref{alpha-i} hold with strict inequality for $i=1,2$,
and the conditions \eqref{a12-1}, \eqref{a12-2} are replaced by the more strict
\be{lin}
\alpha_1 + \alpha_2 > \frac{1}{2} + \frac{n+3}{n-1}(\frac{1}{p}-\frac{1}{2}).
\end{equation}
If $p=p_0$, then one can let \eqref{beta-0} be obeyed with equality.
\end{theorem}

This is sharp except for endpoints when $n=2$ (since \eqref{lin} nearly matches
\eqref{a12-2} in this case), and is somewhat sharp for $n > 2$.  One can replace the multipliers $D_0$, $D_+$, $D_-$ by other symbols which satisfy the same types of regularity and decay estimates, but we shall not pursue this matter here.

\begin{proof} In order to apply the complex interpolation method we shall allow the indices $\beta_0$, $\beta_+$, $\beta_-$ to acquire an imaginary part $it$.  For notational simplicity we keep ourselves to the case $t=0$, but the reader may easily verify that the following argument also works for arbitrary $t$ with constants which are at most exponential in $t$.

When $p=2$ these claims were proven in \cite{damiano:null} (see also \cite{tv:cone2}).  Since we have replaced \eqref{a12-1}, \eqref{a12-2} with the linear condition \eqref{lin}, we see by complex interpolation that it suffices to verify the claims when $p = p_0$.  In this case \eqref{beta-0} becomes $\beta_0 \geq 0$; since the operator $D_0^{\beta_0} D_+^{-\beta_0}$ is bounded in $L^p$ by standard multiplier theory we may assume that $\beta_0 = 0$.  

Let $\phi$, $\psi$ be solutions to the free wave equation.  By a finite decomposition and time reversal symmetry we may assume that $\hat \phi$ is supported in the upper light cone $\Sigma^+$ and that $\hat \psi$ is supported on either the upper or the lower light cone $\Sigma^\pm$.

Write $\phi =: \sum_j \phi_j$, $\psi =: \sum_k \psi_k$, where $\phi_j$, $\psi_k$ have frequency supports on the region $D_+ \sim 2^j$, $D_+ \sim 2^k$ respectively.
We now rewrite \eqref{null} as
$$ \| \sum_{j,k} D_+^{\beta_+ - \beta_-} |\Box|^{\beta_-} (\phi_j \psi_k) \|_p \lesssim (\sum_j 2^{2\alpha_1 j} E(\phi_j))^{1/2}
(\sum_k 2^{2\alpha_2 k} E(\psi_k))^{1/2}.
$$
It suffices to show that
$$ \| D_+^{\beta_+ - \beta_-} |\Box|^{\beta_-} (\phi_j \psi_k) \|_p \lesssim 2^{-\eps |j-k|} (2^{2\alpha_1 j} E(\phi_j))^{1/2}
(2^{2\alpha_2 k} E(\psi_k))^{1/2}
$$
for all $j,k$ and some $\eps > 0$, since the claim then follows from the triangle inequality, Cauchy-Schwarz, and Young's inequality for sums.

By symmetry we may take $j \leq k$; by scaling and \eqref{scaling} we may take $j=0$.  Having used \eqref{scaling}, our arguments will no longer require this condition and we shall discard it.  Since we are assuming \eqref{alpha-i} to hold with strict inequality we may absorb the $\eps$ into the $\alpha_i$ factors, and reduce ourselves to showing that
$$ \| D_+^{\beta_+ - \beta_-} |\Box|^{\beta_-} (\phi_0 \psi_k) \|_p \lesssim 2^{\eps k} 2^{\alpha_2 k} E(\phi_0)^{1/2} E(\psi_k)^{1/2}.
$$

Fix $k$.  For any $l \geq 0$, decompose the double light cone into finitely overlapping projective sectors $\Gamma$ of angular width $2^{-l}$, and let $\phi_0 = \sum_\Gamma \phi_{0,\Gamma}$, $\psi_k = \sum_\Gamma \psi_{k,\Gamma}$ be a Fourier decomposition of $\phi_0$, $\psi_k$ subordinate to these sectors.  It suffices to show that
$$ \| \sum_{\Gamma, \Gamma': \angle(\Gamma,\Gamma') \sim 2^{-l}} 
D_+^{\beta_+ - \beta_-} |\Box|^{\beta_-} (\phi_{0,\Gamma} \psi_{k,\Gamma'}) \|_p \lesssim 2^{-\eps l} 2^{\eps k} 2^{\alpha_2 k} E(\phi_0)^{1/2} E(\psi_k)^{1/2}
$$
for all $l \geq 0$, since the claim follows by summing in $l$ using the bilinear partition of unity based on angular separation (see e.g. \cite{tvv:bilinear}, \cite{tv:cone2}).  By the triangle inequality and Cauchy-Schwarz as before, it suffices to show that
\be{tar} \| 
D_+^{\beta_+ - \beta_-} |\Box|^{\beta_-} (\phi_{0,\Gamma} \psi_{k,\Gamma'}) \|_p \lesssim 2^{-\eps l} 2^{\eps k} 2^{\alpha_2 k} E(\phi_{0,\Gamma})^{1/2} E(\psi_{k,\Gamma'})^{1/2}
\end{equation}
for each $\Gamma$, $\Gamma'$ with angular separation $2^{-l}$.  By a spatial rotation we may assume that $\phi$ has frequency support in \eqref{phi-sec}, and that either $\psi$ or $\overline{\psi}$ has frequency support in \eqref{psi-sec}.

Now suppose $k \gg 1$.  In this case $D_+^{\beta_+-\beta_-}$ is equal to $2^k$ times a harmless multiplier on the frequency support of $\phi_{0,\Gamma} \psi_{k,\Gamma'}$, so by Proposition \ref{toy} the left-hand side of \eqref{tar} is majorized by
$$ 2^{k(\beta_+-\beta_-)} 
2^{\beta_- (k-2l)} 2^{l(\frac{n+1}{p} - (n-1))}
2^{k(\frac{1}{p}-\frac{1}{2}+\eps)} E(\phi_{0,\Gamma})^{1/2} E(\psi_{k,\Gamma'})^{1/2}.$$
The claim \eqref{tar} then follows after some algebra from \eqref{alpha-i}, \eqref{scaling}, and the assumption that \eqref{beta-minus} holds with strict inequality.

It remains to consider the case $k = O(1)$; by a mild Lorentz transformation we may make $k=0$.  (This affects $D_+$ slightly, but this change is irrelevant).
Since $\Gamma$ and $\Gamma'$ differ in angle by $2^{-l}$, some geometry shows that $D_+$ is at least $2^{-l}$ on the frequency support of $\phi \psi$.  We may thus majorize $D_+^{\beta_+ - \beta_-}$ in $L^p$ by $1 + 2^{l (\beta_+ - \beta_-)}$, and
by Proposition \ref{toy} again the left-hand side of \eqref{tar} is bounded by
$$
(1 + 2^{l(\beta_+ - \beta_-)})
2^{\beta_- (-2l)} 2^{l(\frac{n+1}{p} - (n-1))}
E(\phi_{0,\Gamma})^{1/2} E(\psi_{0,\Gamma'})^{1/2}.$$
The claim \eqref{tar} then follows after some algebra from \eqref{lin}, \eqref{p0-def}, \eqref{scaling}, and the assumption that \eqref{beta-minus} holds with strict inequality.
\end{proof}

It seems plausible that many of the missing endpoints in the above result could be obtained if one was willing to prove a generalization of Proposition \ref{toy} which considered the interaction of multiple scales and multiple angular separations using more sophisticated tools than the triangle inequality, and in which the $2^{\eps k}$ loss was eliminated.  

Apart from the issue of endpoints, there is still the unsatisfactory gap between
the condition \eqref{lin} in Theorem \ref{n}, and the necessary conditions \eqref{a12-1}, \eqref{a12-2} conjectured in \cite{damiano:null}, when $n > 2$.
In particular, to resolve the conjecture we would need to consider exponents near the case
\be{endpt}
p = \frac{n+2}{n} = p_0(n-1), \quad \alpha_1 + \alpha_2 = 1/p,
\end{equation}
which is the intersection of \eqref{a12-1} and \eqref{a12-2}.  This case turns out to be related to the open (and quite difficult)

\begin{conjecture}[Machedon-Klainerman for the Schr\"odinger equation]\label{mks} Let $\phi$, $\psi$ be functions on $\R^{n+1}$ which have frequency supports in
$$ \{ (\xi, \frac{1}{2}|\xi|^2): \xi \in \underline \Sigma \}$$
and
$$ \{ (\xi, \frac{1}{2}|\xi|^2): \xi \in -\underline \Sigma \}$$
respectively.  Then \eqref{est-0} holds for all $p \geq p_0(n)$.
\end{conjecture}

\begin{proposition}  Suppose that \eqref{null} holds for some set of exponents satisfying \eqref{endpt}.  Then Conjecture \ref{mks}
holds with $n$ replaced by $n-1$.
\end{proposition}

\begin{proof}
We use the method of descent, exploiting the fact that the paraboloid in $\R^n$ is a conic section of the light cone in $\R^{n+1}$.

To prove Conjecture \ref{mks} for $n-1$, it suffices to verify it for $p = p_0(n-1) = \frac{n+2}{n}$, since the other endpoint $p=\infty$ is trivial.  From Plancherel's theorem it suffices to show that
\be{ta} \| \int_{\underline \Sigma_{n-1}} \int_{\underline \Sigma_{n-1}}
 e^{2\pi i x \cdot (\xi - \eta)}
e^{2\pi i t \cdot \frac{1}{2}(|\xi|^2 + |\eta|^2)}
f(\xi) g(\eta)\ d\xi d\eta \|_{L^p(\R^{n-1 + 1})} \lesssim \|f\|_2 \|g\|_2
\end{equation}
for all $C^\infty$ functions $f,g$ on $\underline \Sigma_{n-1}$.

Fix $f$, $g$.  Let $R\gg 1$ be a large number, and consider the functions $\phi_R$, $\psi_R$ defined on $\R^{n+1}$ by
$$ \phi_R(x_1, x', t) := \int_{-1}^1 \int e^{2\pi i (t+x_1) (R+s)} e^{2\pi i (t-x_1)\frac{|\xi|^2}{4(R+s)}} 
e^{2\pi i x' \cdot \xi} 
f(\xi)\ d\xi ds$$
and
$$ \psi_R(x_1, x', t) := \int_{-1}^1 \int e^{-2\pi i (t+x_1) (R+s)} e^{2\pi i (t-x_1)\frac{|\xi|^2}{4(R+s)}} 
e^{2\pi i x' \cdot \xi} 
g(\xi)\ d\xi ds.$$
One can easily verify that $\phi_R$ and $\psi_R$ solve the wave equation, and that
$$ \| \phi_R[0] \|_{\dot H^{\alpha_1}} \sim R^{\alpha_1} \|f\|_2, \quad
\| \psi_R[0] \|_{\dot H^{\alpha_2}} \sim R^{\alpha_2} \|g\|_2.$$
By the assumption that \eqref{null} holds at \eqref{endpt}, we have
$$ \| D_0^{\beta_0} D_-^{\beta_-} D_+^{\beta_+} (\phi \psi) \|_p
\lesssim R^{1/p} \|f\|_2 \|g\|_2$$
for some $\beta_0$, $\beta_-$, $\beta_+$.  However, a computation of the frequency support of $\phi\psi$ shows that $D_0$, $D_-$, $D_+$ are all given by smooth functions comparable to 1 on this support, and so we have
$$ \| \phi \psi \|_p
\lesssim R^{1/p} \|f\|_2 \|g\|_2$$
Making the change of variables $X = t + x_1$, $T = (t-x_1)/2R$, we have
\bas \phi \psi(X,x',T) = e^{4\pi i XR}
\int_{-1}^1 \int_{-1}^1
\int\int &e^{2\pi i X (s+s')} e^{2\pi i T
(\frac{|\xi|^2}{2 + s/2R} + \frac{|\eta|^2}{2 + s'/2R})} 
e^{2\pi i x' \cdot (\xi + \eta)} \\
&f(\xi) g(\eta)\ d\xi d\eta ds ds'
\end{align*}
and so we have
$$ \| \phi \psi \|_{L^p_X L^p_{x'} L^p_T} \lesssim \|f\|_2 \|g\|_2,$$
the $R^{1/p}$ factor having cancelled against the Jacobian term.  

The phase factor $2^{4\pi i XR}$ can be discarded. Letting $R \to \infty$ and taking limits, and then evaluating the $s$, $s'$ integrations, we thus see that
$$ \| 
(\frac{\sin 2\pi X}{X})^2
\int\int e^{2\pi i T (\frac{|\xi|^2}{2} + \frac{|\eta|^2}{2})} 
e^{2\pi i x' \cdot (\xi + \eta)} 
f(\xi)\ d\xi d\eta \|_{L^p_X L^p_{x',T}} \lesssim \|f\|_2 \|g\|_2,$$
and \eqref{ta} follows since the $X$ behaviour is trivial.
\end{proof}

The Machedon-Klainerman conjecture for the Schr\"odinger equation looks very similar to the results proved in Theorem \ref{main} with $k=0$, however the method of proof breaks down because the strong Huygens' principle (\eqref{phi-eek} and \eqref{redblue-h}) totally fails for this equation.  The arguments in \cite{wolff:cone}
fail for similar reasons (basically, the tubes $T$ are no longer constrained to
point in null directions).  Only partial progress is known for this problem; for instance, when $n=2$ this conjecture should hold for all $p \geq 2 - \frac{1}{3}$, but the best result achieved to date is $p > 2 - \frac{2}{17}$, see \cite{tv:cone1}.
A resolution of this conjecture would also yield new results for the very difficult restriction and Bochner-Riesz problems for the paraboloid \cite{tvv:bilinear}, \cite{carbery:parabola} as well as the problem of pointwise convergence of the Schr\"odinger equation to the initial data \cite{tv:cone2}.

The estimate \eqref{null} at the endpoint \eqref{endpt} shares some features in common with the conjectures (29), (30) in \cite{wolff:cone}, and the resolutions to these problems may well be related.

\section{Appendix I: Proof of Lemma \ref{tube-decomp}}\nonumber

We now give the proof of Lemma \ref{tube-decomp}.  The main technical difficulty is to obtain a good constant in the main term of \eqref{bessel-T}; this will be obtained by using characteristic functions to decompose the Fourier domain, followed by an averaging over rotations to smooth things out spatially.

Partition 
$$ S^{n-1} \cap \underline \Sigma = \bigcup_{\omega \in \Ecal} A_\omega$$
where $A_\omega$ consists of those points in $S^{n-1} \cap \underline \Sigma$ which are closer to $\omega$ than any other element of $\Ecal$.  Thus $A_\omega$ is in the $O(1/r)$-neighbourhood of $\omega$. 

Let $G \subset SO(n)$ denote the set of all rotations in $\R^n$ which differ from the identity by $O(1/r)$.  Let $d\Omega$ be a smooth compactly supported probability measure on the interior of $G$.

For each $\Omega \in G$ and $\omega \in \Ecal$, we define the Fourier projection operators $P_{\Omega,\omega}$ for test functions $f$ on $\R^n$ by
$$ \widehat{P_{\Omega,\omega} f}(\xi) := \chi_{\Omega(A_\omega)}(\xi/|\xi|) 
\hat f(\xi).$$
Note that 
\be{rough-decomp}
\phi(0) = \sum_{\omega \in \Ecal} P_{\Omega,\omega} \phi(0) 
\end{equation}
for all $\Omega \in G$ and red waves $\phi$ of frequency 1. 

The decomposition \eqref{rough-decomp} is well-behaved in frequency, but has terrible spatial localization properties.  To get around this we will now average in $d\Omega$ and then apply a spatial cutoff.

Recall the function $\eta_0$ constructed in Section \ref{localize}.  Write 
$$\eta^{x_0}(x) := \eta_0(\frac{c^2}{r}(x-x_0))$$
for all $x_0 \in L$.  From the Poisson summation formula we have
$$1 = \sum_{x_0 \in L} \eta^{x_0}.$$

We define the functions $\phi_T$ at time 0 by
\be{phit-def} \phi_T(x,0) := \eta^{x_T}(x) \int
(P_{\Omega, \omega_T} \phi(0))(x)\ d\Omega
\end{equation}
and at other times by
$$ \phi_T(t) := U(t) \phi_T(0),$$
where $U(t)$ was the evolution operator defined in \eqref{u-def}.
One may verify that $\phi_T$ is a red wave with angular dispersion $O(1/r)$, that the map $\phi \mapsto \phi_T$ is linear, and that \eqref{margin-phit}, \eqref{recon-1} hold.

We now show \eqref{local2}.  From \eqref{phit-def} and the rapid decay of $\eta^{x_T}$ we have the pointwise estimate
\be{ptwise}
\phi_T(0) \lesssim c^{-C}
\tilde \chi_T(0)^4 \int |P_{\Omega, \omega_T} \phi(0)|\ d\Omega,
\end{equation}
so it suffices by Minkowski and the $L^2$ boundedness of $P_{\Omega,\omega_T}$ to show that
$$ 
\| \tilde \chi_T(0) P_{\Omega, \omega_T} \phi(0) \|_2
\lesssim C_0^C \| \tilde \chi_T(t) P_{\Omega, \omega_T} \phi(t) \|_2.
$$
for all $\Omega \in G$, since the $C_0^C$ factor can be absorbed into the $c^{-C}$ factor by the hypothesis $c \leq 2^{-C_0}$.
In fact, we will show the stronger
\be{schur}
\| \phi \|_{L^2(D(x_D,t_D;r))}
\lesssim C_0 \| \tilde \chi_{D(x_D-(t-t_D)\omega_T,t;r)}^3 \phi(t) \|_2
\end{equation}
for all disks $D$ of radius $r$ and time co-ordinate $t_D = O(R)$, and any red wave $\phi$ with frequency support in the sector $\{ (\xi,|\xi|): |\xi| \sim 1, \angle(\xi,\omega_T) \lesssim C/r \}$.
The previous claim follows by breaking up $\tilde \chi_T(0)$ into various disks $D$.

We now show \eqref{schur}. We have the identity
$$ \phi(x,t_D) = \int  
e^{2\pi i (x \cdot \xi + (t - t_D) |\xi|)} \varphi(\xi) \widehat{\phi(t)}(\xi)\ d\xi$$
where $\varphi$ is an arbitrary bump function which equals one on the sector $\{ \xi \in \underline \Sigma: \angle(\xi,\omega_T) \lesssim 1/r \}$ and is adapted to a slight enlargement of this sector. By standard stationary phase estimates and the observation that $t_D - t = O(C_0 r^2)$, we thus have
$$ \phi(t_D) = \phi(t) * K$$
where the kernel $K$ satisfies the estimates
$$ |K(x)| \lesssim C_0^C r^{1-n} (1 + |x-(t-t_D)\omega_T|/r)^{-N^{10}} (1 + |(x-(t-t_D)\omega_T) \cdot \omega_T|)^{-N^{10}}.$$
The claim \eqref{schur} then follows from a dyadic decomposition of $K$ around the point $(t-t_D) \omega_T$, followed by Young's inequality.  This proves \eqref{local2}.

We now show \eqref{local3}.  Write $R' := \dist(T,Q)$.  Let $D$ denote the disk $D = D(0,0; C^{-1} R')$.  Since $R' \gg R$, it suffices from \eqref{prop}, \eqref{k-bound}, and crude estimates to show that
$$ \| \phi_T \|_{L^2(D)} \lesssim {R'}^{-N} E(\phi)^{1/2}.$$
But this is clear from \eqref{phit-def} and the estimates
$$\| \eta^{x_T} \|_{L^\infty(D)} \lesssim {R'}^{-N}, \quad
\| P_{\Omega,\omega_T} \phi(0)\|_2 \leq E(\phi)^{1/2}.$$
This proves \eqref{local3}.

We now show \eqref{local}.  From \eqref{schur} the left-hand side of \eqref{local} is majorized by
$$
r \sum_T \sup_q \tilde \chi_T(x_q,t_q)^{-3}  
\| \tilde \chi_{D(x_q + t_q \omega_T,0;r)}^{3} \phi_T \|_{2}^2.
$$
From \eqref{ptwise} this is majorized by
$$
c^{-C} r \sum_T \sup_q \tilde \chi_T(x_q,t_q)^{-3}  
\| \tilde \chi_{D(x_q + t_q \omega_T,0;r)}^{3} 
\tilde \chi_T(0)^{4} \int |P_{\Omega, \omega_T} \phi(0)|\ d\Omega
\|_2^2.
$$
From \eqref{tchit-def}, \eqref{tchib-def} we have the elementary inequality
$$ 
\tilde \chi_{D(x_q + t_q \omega_T,0;r)}(x)
\tilde \chi_T(x,0)
\lesssim
\tilde \chi_T(x_q,t_q)
$$
so we may bound the previous by
$$
c^{-C} r \sum_T 
\| 
\tilde \chi_T(0) \int |P_{\Omega, \omega_T} \phi(0)|\ d\Omega
\|_2^2.
$$
Summing in $x_T$, we reduce to showing that
$$
\sum_{\omega \in \Ecal}
\| 
\int |P_{\Omega, \omega} \phi(0)|\ d\Omega
\|_2^2
\lesssim E(\phi).$$
But this follows from Minkowski's inequality followed by Plancherel's theorem.

We now show \eqref{bessel-T}.  We expand the left-hand side as
$$ (\sum_{q_0} \| \int
\sum_{\omega \in \Ecal} \sum_{x_0 \in L}
\m_{q_0,T(\omega,x_0)}
\eta^{x_0}
P_{\Omega,\omega}(\phi(0)) \ d\Omega \|_2^2)^{1/2}.$$
By Minkowski's inequality this is less than or equal to
\be{sook} \int (\sum_{q_0} \| 
\sum_{\omega} \sum_{x_0}
\m_{q_0,T(\omega,x_0)}
\eta^{x_0}
P_{\Omega,\omega}(\phi(0)) \|_2^2)^{1/2}\ d\Omega.
\end{equation}
Define the set $Y \subset S^{n-1}$ by
$$ Y := \bigcup_{\omega \in \Ecal} 
\{ \alpha \in A_\omega: \dist(\alpha, S^{n-1} \backslash A_\omega) > C c^2/r \};$$
this set $Y$ plays a similar role here to the sets $I^{c,k}(Q)$ used in Proposition \ref{lung}.
and define $P_{\Omega(Y)}$ to be the multiplier
$$ \widehat{P_{\Omega(Y)} f}(\xi) := \chi_{\Omega(Y)}(\frac{\xi}{|\xi|}) \hat f(\xi).$$
By the triangle inequality, \eqref{sook} is less than the sum of
\be{sook1} \int (\sum_{q_0} \| 
\sum_{\omega} \sum_{x_0}
\m_{q_0,T(\omega,x_0)}
\eta^{x_0}
P_{\Omega,\omega} P_{\Omega(Y)} \phi(0) \|_2^2)^{1/2}
\ d\Omega.
\end{equation}
and
\be{sook2} \int (\sum_{q_0} \| 
\sum_{\omega} \sum_{x_0}
\m_{q_0,T(\omega,x_0)}
\eta^{x_0}
P_{\Omega,\omega} (1-P_{\Omega(Y)}) \phi(0) \|_2^2)^{1/2}
\ d\Omega.
\end{equation}
Consider the quantity \eqref{sook1}.  The contributions of each $\omega$ are orthogonal as $\omega$ varies, so we can rewrite \eqref{sook1} as
$$ \int (\sum_{q_0} \sum_{\omega}
\| \sum_{x_0} \m_{q_0,T(\omega,x_0)} \eta^{x_0}
P_{\Omega,\omega} P_{\Omega(Y)} \phi(0) \|_2^2)^{1/2}
\ d\Omega,$$
which can be re-arranged as
$$ \int (\sum_\omega \int |(P_{\Omega,\omega} P_{\Omega(Y)} \phi(0))(x)|^2
\sum_{q_0} (\sum_{x_0} \m_{q_0,T(\omega,x_0)} 
\eta^{x_0}(x) )^2\ dx)^{1/2}\ d\Omega.$$
This is clearly less than or equal to
$$ \int (\sum_\omega \int |(P_{\Omega,\omega} P_{\Omega(Y)} \phi(0))(x)|^2
(\sum_{q_0} \sum_{x_0} \m_{q_0,T(\omega,x_0)} 
\eta^{x_0}(x))^2\ dx)^{1/2}\ d\Omega.$$
Summing in $q_0$ and then in $x_0$, this simplifies to
$$ \int (\sum_\omega \int |(P_{\Omega,\omega} P_{\Omega(Y)} \phi(0))(x)|^2\ dx)^{1/2}\ d\Omega.$$
By the orthogonality of the $P_{\Omega,\omega}$, we thus have
$$ \eqref{sook1} \leq \int \| P_{\Omega(Y)} \phi(0) \|_2\ d\Omega \leq E(\phi)^{1/2}.$$

Now consider \eqref{sook2}.  Repeating the above argument, but noting that we must use almost orthogonality instead of orthogonality, we obtain
$$ \eqref{sook2} \lesssim \int \| (1-P_{\Omega(Y)}) \phi(0) \|_2\ d\Omega.$$
By Cauchy-Schwarz we thus have
$$ \eqref{sook2} \lesssim (\int \| (1-P_{\Omega(Y)}) \phi(0) \|_2^2\ d\Omega)^{1/2};$$
by Plancherel we thus have
$$ \eqref{sook2} \lesssim (\int |\hat \phi(\xi,0)|^2 
(\int \chi_{\Omega(Y)^c}(\frac{\xi}{|\xi|})\ d\Omega)\ d\xi)^{1/2}.$$
The inner integral can be seen to be $O(c^2)$ for all $\xi$, so by Plancherel again we have
$$ \eqref{sook2} \lesssim c E(\phi)^{1/2}.$$
Combining our estimates for \eqref{sook1} and \eqref{sook2} we obtain \eqref{bessel-T}.
\endprf

\section{Appendix II: Proof of Lemma \ref{far-cone}}\nonumber

To complete the proof of Theorem \ref{main} we have to prove Lemma \ref{far-cone}.  This shall be a straightforward application of the wave packet decomposition Lemma \ref{tube-decomp} and the energy estimates on cones developed in Section \ref{opposite-sec}.  However, our arguments are not as delicate as those in the rest of the proof of Theorem \ref{main}, as Lemma \ref{far-cone} is not an endpoint estimate.  For instance, we will be able to lose powers of $R^{1/N}$ in our estimates.

We turn to the details.  Let $R \geq 2^{C_0 k}$, and let $\phi$, $\psi$ be red and blue waves of frequency $1$, $2^k$ respectively and margins at least $1/200$.  We may assume the energy normalization \eqref{energy-bound}.  By translation invariance it suffices to show that
$$ \| (P_D \phi) \psi \|_{L^p(Q^{ann}(0,0;R, 2R))},
\| \phi (P_D \psi) \|_{L^p(Q^{ann}(0,0;R, 2R))}
 \lesssim R^{-1/C}$$
for any disk $D$ centered at the origin with radius $2^k \leq r \leq C_0 R^{1/2 + 4/N}$.  We shall only prove the estimate for $(P_D \phi) \psi$; the estimate for $\phi (P_D \psi)$ is proven similarly, observing that any additional powers of $2^k$ which appear can be absorbed in to the $R^{-1/C}$ factor.

As in previous arguments, we obtain the $L^p$ estimate by interpolation between an $L^1$ estimate and an $L^2$ estimate, namely
\be{l1-cloc}
\| (P_D \phi) \psi \|_{L^1(Q^{ann}(0,0;R, 2R))},
\lesssim R^{C/N} R^{3/4}
\end{equation}
and
\be{l2-cloc}
\| (P_D \phi) \psi \|_{L^2(Q^{ann}(0,0;R, 2R))},
\lesssim R^{C/N} R^{-(n-1)/4}.
\end{equation}
By interpolation and some algebra involving \eqref{p-def} one obtains the desired result if $N$ is sufficiently small.  The point is that \eqref{l1-cloc} improves over \eqref{l1-noloc} by a substantial power of $R$.

We first prove \eqref{l1-cloc}.  The contribution outside of the set $C^{red}(0,0;R^{1/2})$ is acceptable by \eqref{phi-eek}, and the contribution inside this set is acceptable by Corollary \ref{doublecone}.  This proves \eqref{l1-cloc}.

Now we prove \eqref{l2-cloc}.  Set $Q := Q(0,0;2R)$.  We use Lemma \ref{tube-decomp} with this cube $Q$ and some arbitrary value of $c$ (say $c = 2^{-C_0}$) to decompose $(P_D \phi)$ as
$$ (P_D \phi) = \sum_{T \in \Tcal} (P_D \phi)_T.$$

By \eqref{local3} and crude estimates (e.g. \eqref{l2-noloc} or Lemma \ref{mock}) the contribution of those tubes $T$ for which $\dist(T, 0) \gg R$ is acceptable, so we can restrict ourselves to those tubes $T$ for which $\dist(T,0) \lesssim R$.  There are only $O(R^C)$ of these tubes.  

We now dispose of those remaining tubes $T$ for which $R^{1/2 + 1/N} \lesssim \dist(T, 0) \lesssim R$.  For such tubes $T$ we see from \eqref{local2} with $t=0$ and \eqref{essential-concentration} that 
$$ E((P_D \phi)_T) \lesssim \| \tilde \chi_T(0) P_D \phi(0) \|_2^2 \lesssim
R^{-N} E(\phi) = R^{-N},$$
and one can dispose of the contribution of these tubes by crude estimates.

It remains to deal with those tubes for which $\dist(T,0) \lesssim R^{1/2 + 1/N}$.  In this case we want to exploit the fact that $(P_D \phi)_T$ is concentrated on $R^{1/N} T$, the dilate of $T$ around its center by $R^{1/N}$.  Indeed, from \eqref{local} we see that
$$ \| (P_D \phi)_T \|_{L^2(Q \backslash R^{1/N} T)} \lesssim R^{-N},$$
and so the contribution outside $R^{1/N} T$ is acceptable by crude estimates.  It thus remains to control the expression
$$ \| \sum_T (P_D \phi)_T \psi \chi_{R^{1/N} T} \|_{L^2(Q^{ann}(0,0;R, 2R))}.$$
We now observe the geometric fact that as $T$ ranges over all tubes in $\T$ with $\dist(T,0) \lesssim R^{1/2 + 1/N}$, the tubes $\chi_{R^{1/N} T}$ have an overlap of at most $O(R^{C/N})$ in $Q^{ann}(0,0;R, 2R))$.  By almost orthogonality, it thus suffices to show that
$$ (\sum_T \| (P_D \phi)_T \psi \|_2^2)^{1/2} \lesssim R^{-(n-1)/4}.$$
Applying Lemma \ref{mock}, we can bound the left-hand side by
$$ R^{-(n-1)/4} (\sum_T E((P_D \phi)_T) E(\psi))^{1/2}.$$
But this is acceptable by \eqref{bessel} and \eqref{energy-minor}.
\endprf

\end{document}